\documentclass[a4paper,1p,sort&compress,dvips]{elsarticle}

\usepackage{stmaryrd}
\usepackage{latexsym}
\usepackage{amsmath}
\usepackage{amsfonts}
\usepackage{amssymb}
\usepackage{bm}
\usepackage[english]{babel}

\DeclareMathAlphabet{\mathpzc}{OT1}{pzc}{m}{it}

\journal{{\tt arXiv}, typeset with elsarticle.cls}
\newcommand{\be}{\begin{equation}}
\newcommand{\ee}{\end{equation}}

\newcommand{\ef}[1]{\, #1}

\newcommand{\gc}{{\cal G}}

\newcommand{\ssm}{\smallsetminus}
\renewcommand{\emptyset}{\varnothing}

\newcommand{\fpl}{\mathpzc{Fpl}}
\newcommand{\LP}{\mathpzc{LP}}

\newcommand{\Sym}{\mathrm{Sym}}
\newcommand{\RS}{\mathsf{RS}}

\newcommand{\sqqV}{\,\rule{.5pt}{7pt}\rule{6pt}{0pt}\rule{.5pt}{7pt}\,}
\newcommand{\sqqH}{\,\makebox[0pt][l]{\rule{7pt}{.5pt}}\raisebox{6.5pt}{\rule{7pt}{.5pt}}\,}

\newcommand{\nn}{\mathcal{N}}

\newcommand{\kket}[1]{\mbox{$\| #1 \rangle \! \rangle $}}

\newcommand{\ket}[1]{\mbox{$| #1 \rangle $}}

\newcommand{\nb}{n_{\bullet}}
\newcommand{\nc}{n_{\circ}}

\newcommand{\ssc}[2]{\ket{s^{#1}_{\circ, #2}}}
\newcommand{\ssb}[2]{\ket{s^{#1}_{\bullet, #2}}}

\newcommand{\ssbcx}[1]{\ket{x_{#1}}}

\newtheorem{corollary}{Corollary}[section]
\newtheorem{conj}{Conjecture}[section]
\newtheorem{lemma}{Lemma}[section]

\newtheorem{proposition}{Proposition}[section]
\newtheorem{definition}{Definition}[section]

\begin{document}

\begin{frontmatter}

\title{Proof of the Razumov-Stroganov conjecture}

\author[ens]{Luigi Cantini}
\ead{lui{}gi.ca{}nti{}ni@lpt.e{}ns.fr}
\author[mi]{Andrea Sportiello}
\ead{Andr{}ea.Sportie{}llo@mi.in{}fn.it}

\address[ens]{Laboratoire de Physique Th\'eorique, \'Ecole Normale
  Sup\'erieure,\\
  24 rue Lhomond, 75005 Paris, France}
\address[mi]{Dipartimento di Fisica dell'Universit\`a degli Studi di
  Milano, and INFN,\\
  via Celoria 16, 20133 Milano, Italy}

\date{\today}

\begin{abstract}
The Razumov-Stroganov conjecture relates the ground-state coefficients
in the even-length dense $O(1)$ loop model
to the enumeration of fully-packed loop configuration on the square,
with alternating boundary conditions, refined according to the link
pattern for the boundary points.

Here we prove this conjecture, by mean of purely combinatorial
methods. The main ingredient is a generalization of the Wieland proof
technique for the dihedral symmetry of these classes, based on the
`gyration' operation, whose full strength we will investigate in a
companion paper.
\end{abstract}


\end{frontmatter}

\section{Introduction}

The study of Alternating Sign Matrices (ASM), i.e.\ matrices with
entries $1$, $-1$ and $0$ such that each row and column sums to $1$,
and $1$ and $-1$ alternate along rows and columns, has a long
tradition.  These objects were introduced by Mills, Robbins and Rumsey
\cite{MRR82,MRR83}, motivated by the study of
\emph{$\lambda$-determinants}. The authors recognized immediately the
relation of the resulting enumeration with the ones of several other
problems, most notably \emph{Plane Partitions}, i.e.\ rhombus tilings
of portions of the triangular lattice.

The first proof of their enumeration has been given by Zeilberger
\cite{zeil}, in a sort of \emph{tour de force} by which he essentially
proved, through non-bijective techniques of generating functions, that
ASMs are equinumerous to \emph{Totally Symmetric Self-Complementary
  Plane Partitions} (TSSCPP), whose enumeration formula was previously
proven by Andrews \cite{andrews}.  Slightly later, Kuperberg
\cite{kupe} found a simpler proof which exploited the bijection
between ASM and configurations of the \emph{six-vertex model} with
domain wall boundary condition, a Yang-Baxter integrable system in
statistical mechanics \cite{baxter}.  It was the integrability of this
latter that allowed physicists to come out with an explicit
determinantal formula for its partition function \cite{IK_I}, which
was used in Kuperberg proof. It is worth mentioning that, although not
used in this proof, the specialization of the six-vertex model
pertinent to the uniform measure over all ASM leads to an even
stronger symmetry, and a formula for the partition function that
involves a \emph{Schur function}, for a certain ``triangular'' Young
diagram \cite{StroSchur, Okada}.

Another incarnation of the ASM are the \emph{fully-packed loop}
configurations (FPL) on regions of the square lattice. A FPL is a
colouring, in two colours (say, black and white), of the edges of the
domain,
such that each vertex is adjacent to two edges of each colour. When
the region is a square, and the colouring of the edges of the boundary
is fixed in an alternating fashion, then the FPL are in bijection with
ASM. The reformulation of ASM in terms of FPL leads naturally to
consider enumerations of family of ASM, whose lines of given colour,
in the FPL formulation, present a given connectivity pattern (called
\emph{link pattern}).  The first striking property of these
enumerations, noted by Bosley and Fidkowski and proven by Wieland
\cite{wie}, is that they are symmetric under a dihedral symmetry
$D_{2n}$ (for a square of side $n$), much larger than the obvious
symmetry group for FPL on the square.

A much stronger fact was pushed forward by Razumov and Stroganov
\cite{rs}, who conjectured that the the enumerations of FPL with a
given link pattern appear as components of the \emph{ground-state
  wavefunction} in the dense $O(1)$ loop model on a semi-infinite
cylinder (a \emph{different} Yang-Baxter integrable model), i.e., the
steady state w.r.t.\ the Markov Chain associated to the transfer
matrix of the model.  Besides the striking numerical evidence in
favour of the conjecture, several particular cases have been solved
positively in the literature. Among these, the sum rule was proven by
Di Francesco and Zinn-Justin \cite{pdf-pzj1}, and, for some infinite
families of link patterns it is possible to compare explicit formulae
for FPL enumerations \cite{pdf-pzj-jbz} with exact results on the
$O(1)$ loop model side \cite{pzj-family}.  More generally, up to now,
promising research lines for proving the conjecture have been mainly
lying on the attempt of ``computing'' the FPL enumerations, and
comparing the result with the components of the loop model ground
state \cite{thapper, zj-triangle}, a strategy that, interestingly, has
seen the emergence of the combinatorics of \emph{Littlewood-Richardson
  coefficients} \cite{PZJ_LR, nadeau}.

In the present paper we give a purely combinatorial proof of the
Razumov-Stroganov conjecture.  The main idea is to recognize the
fundamental role of \emph{gyration}, an operation that can be
performed on FPL, which was already introduced by Mills, Robbins and
Rumsey \cite{MRR83} and was the key in Wieland's proof of the larger
dihedral symmetry \cite{wie}.  

A more striking evidence of the role of gyration is in a fact that we
noticed \emph{before} performing the present work, and plan to
illustrate in a longer companion paper \cite{usprep}: the
Razumov-Stroganov conjecture remains true, apart for a global
multiplicative factor, on a large family of more general domains, as
long as these domains are such that the gyration operation induces
dihedral symmetry (cfr.\ figure \ref{fig.multicone}, left, in Section
\ref{sec.persp}, for an illustration).  As a result, we have a
\emph{family} of Razumov-Stroganov conjectures, indicized by various
other integer parameters, besides the size parameter $n$. This raised
the quest for an unified understanding of the conjecture, on this
whole family of domains. As gyration was the tool for classifying
the family, we expected (and it happened to be the case) that it would
have also played a major role in the unified simultaneous
proof~\cite{usprep}.

In \cite{usprep} we will also deal with the case of FPL with
symmetries, for which there exist variants of the Razumov-Stroganov
conjecture \cite{Raz-Str-2}.  This point is discussed more extensively
in a conclusive section, sec.~\ref{sec.persp}.


The paper is organized as follows. In section \ref{sec.statem} we give
precise definitions of the combinatorial objects we deal with. We
introduce the Temperley-Lieb algebra acting on link patterns, and we
formulate the Razumov-Stroganov conjecture. In section \ref{sec.proof}
we show that the conjecture is a consequence of another striking
enumeration symmetry of FPL (to our knowledge previously unnoticed),
Lemma \ref{thm.HVgen} (an illustration of this fact is in figure
\ref{fig.cor.HV}), and a proposition (Prop.~\ref{prop.equiv}) on how
the Razumov-Stroganov conjecture can be reduced to a special case of
this lemma.  Proposition \ref{prop.equiv} will take us some work to
be proven. This is done, \emph{assuming} certain ``gyration
relations'', in Section \ref{sec.proof}. The gyration relations are
proven separately in Section \ref{sec.wie}. Indeed, they come out as a
very special corollary of a broader analysis of gyration, performed in
Section \ref{sec.wie} in a somewhat larger generality w.r.t.\ what
would suffice for the required gyration relations, and will be
performed in an even larger generality in~\cite{usprep}.

The reader may find useful a glossary of definitions reported in
\ref{app.gloss}.

\section{Statement of the conjecture}
\label{sec.statem}

\subsection{Fully-packed loops on the square lattice}
\label{sec.defFPL}

\noindent
Consider a region $\Lambda$ of the square lattice, determined through
a closed path on the dual lattice. This identifies a set of internal
vertices and edges, $V(\Lambda)$ and $E_0(\Lambda)$, and a set of
``boundary'' edges $E_1(\Lambda)$. Call $E=E_0 \cup E_1$ and $2 N$ the
cardinality of $E_1$ (every closed path on the square lattice has even
length).

We are interested in ensembles of configurations $\phi: E \to
\{b,w\}^E$ (black and white) of edge-colourations, satisfying the
\emph{ice rule}: each vertex $v\in V$ is adjacent to two black and two
white edges. We call such a configuration a \emph{fully-packed loop
  configuration} (FPL), and denote with $\fpl(\Lambda)$ this
ensemble. Consider the partition of $\fpl(\Lambda)$ into sub-ensembles
accordingly to the boundary conditions $\tau$ for $\phi$, encoded as
vectors in $\{b,w\}^{E_1}$. We denote by $\fpl(\Lambda; \tau)$ the
ensemble of FPL $\phi$ whose restriction to $E_1$ is $\tau$.

\begin{figure}[tb!]
\begin{center}
\setlength{\unitlength}{15pt}
\begin{picture}(9,8)
\put(0,0){\includegraphics[scale=3]{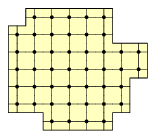}}
\end{picture}
\quad
\begin{picture}(9,8)
\put(0,0){\includegraphics[scale=3]{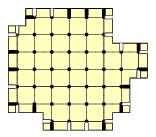}}
\end{picture}
\\
\begin{picture}(9,8)
\put(0,0){\includegraphics[scale=3]{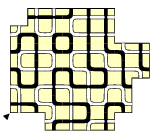}}
\end{picture}
\quad
\begin{picture}(9,8)
\put(0.75,0){\includegraphics[scale=1.25]{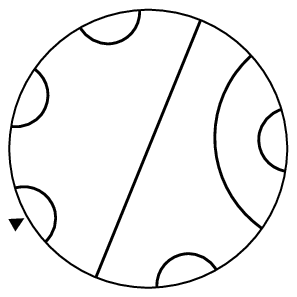}}
\end{picture}
\end{center}
\caption{Top left: an example of domain $\Lambda$. Top right: an
  example of domain with given boundary condition, $(\Lambda,
  \tau)$. Bottom left: an example of FPL $\phi \in \fpl(\Lambda,
  \tau)$. Bottom right: the associated link pattern $\pi(\phi) \in
  \LP(n)$; a small arrow matches the indicization of the endpoints on
  the FPL and on the link pattern.
\label{fig.fplEx}}
\end{figure}

A given $\tau$ has certain sets $E_b(\tau), E_w(\tau) \subseteq E_1$
of black and white entries.  It is easily seen that, if their
cardinalities are odd, then $\fpl(\Lambda; \tau) = \emptyset$. So we
can write $|E_b(\tau)|=2n$ and $|E_w(\tau)|=2(N-n)$.

Because of the ice rule, a configuration $\phi \in \fpl(\Lambda;
\tau)$ causes the set $E(\Lambda)$ to decompose into black and white
closed cycles, and black and white open paths, with endpoints
respectively in $E_b$ and $E_w$. Black paths among themselves, and
white paths among themselves, are non-crossing, while black and white
paths may cross with each other. Label with indices from $1$ to $2n$
the points of $E_b$, in cyclic order.  To a certain FPL $\phi$ we can
thus associate a pairing $\pi(\phi) \in \LP(n)$ of the endpoints,
where $\LP(n)$ is the set of \emph{link patterns}, i.e.\ non-crossing
matchings on the disk, for $2n$ points on the border.  The pairing is
non-crossing, as the square lattice is planar, and the endpoints
are on the boundary of the domain.  We call $\Psi_{\Lambda;
  \tau}(\pi)$ the number of configurations in $\fpl(\Lambda; \tau)$
with link pattern $\pi$.  

Remark that, in order to be definite in the description of
$\pi(\phi)$, we have to specify, besides $\Lambda$, also a cyclic
labeling for the black terminations. Even if we agree on using
counter-clockwise labeling, we have to specify a starting point. We
will be careful on this aspect, all along the paper and within its
figures.

An example of FPL is shown in figure~\ref{fig.fplEx}.

As we said, a simple bijection relates FPL configurations to
configurations in the statistical ensemble of the six-vertex model.
The jargon of this model
suggests to denote by the letters $a$, $b$ and $c$ the six
possible configurations of $\phi$ in a neighbourhood of a vertex,
according to the following rule
(cfr.\ e.g.\ \cite[pp.33-34]{PZJthesis})
\be
\label{eq.defabc}
\setlength{\unitlength}{20pt}
\raisebox{-4mm}{
\begin{picture}(13.5,1.4)(0,-0.2)
\put(0,0){\includegraphics[scale=4]{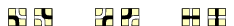}}
\put(0,0.4){$a:$}
\put(5,0.4){$b:$}
\put(10,0.4){$c:$}
\end{picture}
}
\ee
For sets $A \subseteq B$, and $x \in \{b,w\}^B$,
$x|_A$ denotes the restriction of $x$ to the space $\{b,w\}^A$.
For example, $\phi \in \fpl(\Lambda,\tau)$ iff $\phi|_{E'} = \tau$.
Also, for vectors $x \in \{b,w\}^A$, a bar denotes the
\emph{complementation} involution $b \leftrightarrow w$,
i.e.\ $\bar{x}$ is the vector such that $\bar{x}_i = b \leftrightarrow
x_i=w$ and $\bar{x}_i = w \leftrightarrow x_i=b$.
For example, if $\phi \in \fpl(\Lambda, \tau)$, then $\bar{\phi} \in
\fpl(\Lambda, \bar{\tau})$.

\begin{figure}[tb!]
\begin{align*}
&
\setlength{\unitlength}{15pt}
\begin{picture}(7,6)(0,0.5)
\put(0,0){\includegraphics[scale=3]{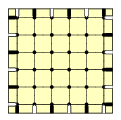}}
\end{picture}
&&
\setlength{\unitlength}{15pt}
\begin{picture}(7,6)(0,0.5)
\put(0,0){\includegraphics[scale=3]{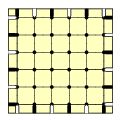}}
\end{picture}
\end{align*}
\caption{Left: the square domain with boundary conditions
  $\tau_+$. Right: the one with boundary
  conditions~$\tau_-$.\label{fig.fplBC}}
\end{figure}

A specially interesting case of domain is the one in which $\Lambda$
is a square of side $n$, and $\tau=\tau_+=(bwbw\ldots bw)$, or the
complementary choice $\tau_-=\overline{\tau_+}=(wbwb\ldots wb)$ (remark that there is no
collision of notation here with $|E_b|=2n$).  The corresponding
domains are shown in figure \ref{fig.fplBC}.  A complete discussion of
this situation, in the framework of interest for this work, can be
found in \cite{DeGier02, PZJthesis}.  We denote by $\fpl(n,\pm)$ the
corresponding ensembles, and $\Psi_{n; \pm}(\pi)$ the corresponding
cardinalities of the refined classes.

In this case, 
a bijection exists with Alternating Sign Matrices
\cite{MRR82, MRR83, bress, faces}, and remarkable combinatorial
relations arise, some of which are proven, others having striking
numerical evidence.  Some examples are
\begin{itemize}
\item A large dihedral symmetry (proven in \cite{wie}), stating that
  $\Psi_{n; \pm}(\pi)$ is invariant under cyclic permutations acting
  over $\pi$, and also that
\be
\label{eq.wiePM}
\Psi_{n; +}(\pi) = \Psi_{n; -}(\pi)
\ee
  (we will thus drop the index $\pm$ in the following).
\item Round formulas for some of these enumerations, among which,
  notably, the cardinality of the whole set,
$
A_n = |\fpl(n)| = \prod_{j=0}^{n-1} \frac{(3j+1)!}{(n+j)!}
$
  (conjectured in \cite{MRR82} and proven in \cite{zeil,kupe}).
\item Identities for special configurations, among which
  $\Psi_{n}(\pi)=A_{n-1}$ for the link pattern
  $\pi=\big((12),(34),(56), \ldots \big)$ (conjectured by J.~Propp in
  \cite{faces}), and various others in \cite{faces, wie}.
\item Polynomiality in $k$ of quantities 
  $\Psi_{n+2k}\big( \pi \!\!\fatslash_k \!\fatbslash_k
\big)$, where 
  $\pi \!\!\fatslash_k \!\fatbslash_k
$ denotes a link pattern $\pi \in \LP(n)$,
  adjoined of a ``rainbow'' of arcs connecting $2n+i$ with
  $2n+2k+1-i$, for $i=1, \ldots, k$ (conjectured in \cite{zubersomenew}
  and proven in \cite{kratcas, kratcasBN}).
\end{itemize}

\noindent
This is the framework of the Razumov-Stroganov conjecture
\cite{rs}. More precisely, the conjecture states the identity (up to a
single normalization overall) between the refined enumerations
$\Psi_n(\pi)$, and a certain set of integers $\tilde{\Psi}_n(\pi)$
arising as components of the ground state of the dense $O(1)$ loop
model, for a cylindric geometry with $2n$ sites per row.

This is a problem arising in the physics of integrable quantum
one-dimensional systems, which started from the context of the XXZ
Quantum Spin Chain, at anisotropy parameter $\Delta=-1/2$, and it
would take us a long detour to give here an approriate introduction
(we refer the reader to \cite{DeGier02, PZJthesis}). Nonetheless, it
is relatively easy to give a purely combinatorial formulation of the
``dense-loop model side'' of the conjecture, at the only price of
introducing a simple diagram algebra acting on the space
$\LP(n)$. This algebra is a representation of the ``affine
Temperley-Lieb Algebra over $2n$ generators, with parameter
$q=e^{\frac{2i\pi}{3}}$ (i.e.\ at a cubic root of unity)'', and, with
some sloppiness, we just call it \emph{Temperley-Lieb Algebra} in the
present context.


\subsection{Temperley-Lieb Algebra}
\label{sec.TL}

\noindent
For $\pi$ a link pattern in $\LP(n)$, and $1 \leq i \leq 2n$, define
$\pi(i)$ as the index matched to $i$. Use cyclic notation for the
indices ($i \equiv i+2n$).

Call $R$ the operator that rotates a link pattern $\pi$ one step
counter-clockwise, or, equivalently, keeps $\pi$ fixed and rotates the
labels one step clockwise
\be
\begin{split}
\pi 
&=
\big( (i_1, j_1), (i_2, j_2), \ldots \big)
\quad \longleftrightarrow \quad
R \pi 
=
\big( (i_1-1, j_1-1), (i_2-1, j_2-1), \ldots \big)
\ef.
\end{split}
\ee
Clearly, $R^{2n}=1$, and $R$ is invertible. Define the $2n$ maps 
$\{ e_j \}_{1 \leq j \leq 2n}$ acting over $\LP(n)$:
\be
\label{eq.78665465734}
e_j (\pi) =
\left\{
\begin{array}{ll}
\pi & \pi(j)=j+1; \\
\!
\displaystyle{
\genfrac{}{}{0pt}{}{\pi \ssm \{ (j,\pi(j)), (j+1,\pi(j+1)) \} \qquad}
{\qquad    \cup \; \{ (j,j+1), (\pi(j),\pi(j+1)) \}}}
    & \mathrm{otherwise}
\end{array}
\right.
\ee
In words, $e_j$ does nothing on $\pi$ if $(j,j+1) \in \pi$, otherwise
it connects $j$ to $j+1$, and the indices previously matched to $j$
and $j+1$ with each other.

These operators are easily seen to
satisfy the following rules
\begin{subequations}
\label{eqs.TL}
\begin{align}
\label{eq.34657865}
e_i
&= R e_{i+1} R^{-1}
\ef;
\\
\label{eq.762347658}
e_i^2
&=
e_i
\ef;
\\
e_i e_j 
&= 
e_j e_i
&&
|i-j| \neq 1
\\
e_i e_{i \pm 1} e_i 
&= e_i
\ef.
\end{align}
\end{subequations}
These rules are deduced by recognizing that, if the link patterns in
$\LP(n)$ are graphically represented as, e.g.,
\[
\pi=
\big( (1,6),(2,3),(4,5),(7,10),(8,9) \big)
\ :
\quad
\setlength{\unitlength}{10pt}
\raisebox{-10pt}{
\begin{picture}(11,4)
\put(0,1){\includegraphics[scale=2]{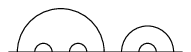}}
\put(1,0){\makebox[0pt][c]{$\scriptstyle{1}$}}
\put(2,0){\makebox[0pt][c]{$\scriptstyle{2}$}}
\put(3,0){\makebox[0pt][c]{$\scriptstyle{3}$}}
\put(4,0){\makebox[0pt][c]{$\scriptstyle{4}$}}
\put(5,0){\makebox[0pt][c]{$\scriptstyle{5}$}}
\put(6,0){\makebox[0pt][c]{$\scriptstyle{6}$}}
\put(7,0){\makebox[0pt][c]{$\scriptstyle{7}$}}
\put(8,0){\makebox[0pt][c]{$\scriptstyle{8}$}}
\put(9,0){\makebox[0pt][c]{$\scriptstyle{9}$}}
\put(10,0){\makebox[0pt][c]{$\scriptstyle{10}$}}
\end{picture}
}
\]
then the action of $R$
and $e_j$ over $\LP(n)$ is graphically encoded by the diagrams
\begin{align}
R:&
\quad
\setlength{\unitlength}{10pt}
\raisebox{-10pt}{
\begin{picture}(11,2.5)
\put(0,1){\includegraphics[scale=2]{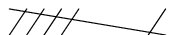}}
\put(1,0){\makebox[0pt][c]{$\scriptstyle{1}$}}
\put(2,0){\makebox[0pt][c]{$\scriptstyle{2}$}}
\put(3,0){\makebox[0pt][c]{$\scriptstyle{3}$}}
\put(4,0){\makebox[0pt][c]{$\scriptstyle{\cdots}$}}
\put(10,0){\makebox[0pt][c]{$\scriptstyle{2n}$}}
\put(7,1.8){\makebox[0pt][c]{$\scriptstyle{\cdots}$}}
\end{picture}
}
\\
e_j:&
\quad
\setlength{\unitlength}{10pt}
\raisebox{-10pt}{
\begin{picture}(11,3.5)
\put(0,1){\includegraphics[scale=2]{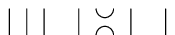}}
\put(1,0){\makebox[0pt][c]{$\scriptstyle{1}$}}
\put(2,0){\makebox[0pt][c]{$\scriptstyle{2}$}}
\put(3,0){\makebox[0pt][c]{$\scriptstyle{3}$}}
\put(4,0){\makebox[0pt][c]{$\scriptstyle{\cdots}$}}
\put(5.7,0){\makebox[0pt][c]{$\scriptstyle{j}$}}
\put(7.2,0){\makebox[0pt][c]{$\scriptstyle{j+1}$}}
\put(9,0){\makebox[0pt][c]{$\scriptstyle{\cdots}$}}
\put(10,0){\makebox[0pt][c]{$\scriptstyle{2n}$}}
\put(4,1.5){\makebox[0pt][c]{$\scriptstyle{\cdots}$}}
\put(9,1.5){\makebox[0pt][c]{$\scriptstyle{\cdots}$}}
\end{picture}
}
\end{align}
The two cases of equation (\ref{eq.78665465734}) are well illustrated
by the action of $e_1$ and $e_2$ over the link pattern above, that
give
\begin{align}
e_1 \pi
\ &:
\quad
\setlength{\unitlength}{10pt}
\raisebox{-10pt}{
\begin{picture}(11,5.5)
\put(0,2.5){\includegraphics[scale=2]{fig_LPopen1.eps}}
\put(0,1){\includegraphics[scale=2]{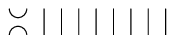}}
\put(1,0){\makebox[0pt][c]{$\scriptstyle{1}$}}
\put(2,0){\makebox[0pt][c]{$\scriptstyle{2}$}}
\put(3,0){\makebox[0pt][c]{$\scriptstyle{3}$}}
\put(4,0){\makebox[0pt][c]{$\scriptstyle{4}$}}
\put(5,0){\makebox[0pt][c]{$\scriptstyle{5}$}}
\put(6,0){\makebox[0pt][c]{$\scriptstyle{6}$}}
\put(7,0){\makebox[0pt][c]{$\scriptstyle{7}$}}
\put(8,0){\makebox[0pt][c]{$\scriptstyle{8}$}}
\put(9,0){\makebox[0pt][c]{$\scriptstyle{9}$}}
\put(10,0){\makebox[0pt][c]{$\scriptstyle{10}$}}
\end{picture}
}
=
\raisebox{-10pt}{
\begin{picture}(11,4)
\put(0,1){\includegraphics[scale=2]{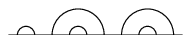}}
\put(1,0){\makebox[0pt][c]{$\scriptstyle{1}$}}
\put(2,0){\makebox[0pt][c]{$\scriptstyle{2}$}}
\put(3,0){\makebox[0pt][c]{$\scriptstyle{3}$}}
\put(4,0){\makebox[0pt][c]{$\scriptstyle{4}$}}
\put(5,0){\makebox[0pt][c]{$\scriptstyle{5}$}}
\put(6,0){\makebox[0pt][c]{$\scriptstyle{6}$}}
\put(7,0){\makebox[0pt][c]{$\scriptstyle{7}$}}
\put(8,0){\makebox[0pt][c]{$\scriptstyle{8}$}}
\put(9,0){\makebox[0pt][c]{$\scriptstyle{9}$}}
\put(10,0){\makebox[0pt][c]{$\scriptstyle{10}$}}
\end{picture}
}
\\
\label{eq.6576488765}
e_2 \pi
\ &:
\quad
\setlength{\unitlength}{10pt}
\raisebox{-10pt}{
\begin{picture}(11,5.5)
\put(0,2.5){\includegraphics[scale=2]{fig_LPopen1.eps}}
\put(0,1){\includegraphics[scale=2]{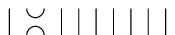}}
\put(1,0){\makebox[0pt][c]{$\scriptstyle{1}$}}
\put(2,0){\makebox[0pt][c]{$\scriptstyle{2}$}}
\put(3,0){\makebox[0pt][c]{$\scriptstyle{3}$}}
\put(4,0){\makebox[0pt][c]{$\scriptstyle{4}$}}
\put(5,0){\makebox[0pt][c]{$\scriptstyle{5}$}}
\put(6,0){\makebox[0pt][c]{$\scriptstyle{6}$}}
\put(7,0){\makebox[0pt][c]{$\scriptstyle{7}$}}
\put(8,0){\makebox[0pt][c]{$\scriptstyle{8}$}}
\put(9,0){\makebox[0pt][c]{$\scriptstyle{9}$}}
\put(10,0){\makebox[0pt][c]{$\scriptstyle{10}$}}
\end{picture}
}
=
\raisebox{-10pt}{
\begin{picture}(11,4)
\put(0,1){\includegraphics[scale=2]{fig_LPopen1.eps}}
\put(1,0){\makebox[0pt][c]{$\scriptstyle{1}$}}
\put(2,0){\makebox[0pt][c]{$\scriptstyle{2}$}}
\put(3,0){\makebox[0pt][c]{$\scriptstyle{3}$}}
\put(4,0){\makebox[0pt][c]{$\scriptstyle{4}$}}
\put(5,0){\makebox[0pt][c]{$\scriptstyle{5}$}}
\put(6,0){\makebox[0pt][c]{$\scriptstyle{6}$}}
\put(7,0){\makebox[0pt][c]{$\scriptstyle{7}$}}
\put(8,0){\makebox[0pt][c]{$\scriptstyle{8}$}}
\put(9,0){\makebox[0pt][c]{$\scriptstyle{9}$}}
\put(10,0){\makebox[0pt][c]{$\scriptstyle{10}$}}
\end{picture}
}
\end{align}
The parameter $q$, here set to
$e^{\frac{2i\pi}{3}}$, would have appeared in the first case of
(\ref{eq.78665465734})
(if $\pi(j) = \hbox{$j+1$}$, then $e_j (\pi)$ also produces an
overall factor $-q-q^{-1}$) and in equation
(\ref{eq.762347658}) (which would read
$e_i^2 =\left(-q-q^{-1} \right) e_i$).
In a graphical representation, as in the pictures of equation
(\ref{eq.6576488765}), we can think to this factor as associated to
the cycles that are detached from the boundary by the diagram action
of $e_j$.  Of course, $-q-q^{-1}=1$ for $q = e^{\frac{2i\pi}{3}}$.

\subsection{A remark on vector notation}
\label{ssec.VectNot}

\noindent
We will adopt all along the paper a ``vector'' notation. Indeed,
various facts we deal with here take the form
\be
\label{eq.7675457}
\forall \  \pi \in \LP(n)
\qquad
A(\pi) = B(\pi)
\ee
for $A(\pi)$ and $B(\pi)$ ``numbers'' associated to
link-pattern configurations~$\pi$.

Such a statement can be phrased in terms of formal vectors
$\ket{\pi}$, taken as the canonical basis of a linear space over the
field $\mathbb{C}$ (or any other field in which $A(\pi)$ and $B(\pi)$
are valued, such as $\mathbb{R}$ or $\mathbb{Q}$).  The dimension of
this linear space is $|\LP(n)|=C_n$, the $n$-th Catalan number, and we
will denote the space as $\mathbb{C}^{\LP(n)}$.  Calling ${\bf 0}$ the
zero vector in this space, the relation above reads
\be
\sum_{\pi \in \LP(n)}
\big( A(\pi) - B(\pi) \big) \ket{\pi}
=
{\bf 0}
\ef.
\ee
If it is understood that 
$\ket{A} = \sum_{\pi} A(\pi) \ket{\pi}$
and
$\ket{B} = \sum_{\pi} B(\pi) \ket{\pi}$, then the identity (\ref{eq.7675457}) is
just the fact that $\ket{A} = \ket{B}$ as vectors in 
this space.

An example of this notation is the statement of the dihedral symmetry.
If $R$ is the rotation operator, we have
\be
\label{eq.763765876}
\forall \  \pi \in \LP(n)
\qquad
\Psi_n(\pi) = \Psi_n(R \pi)
\ee
and can be rephrased as
\be
\sum_{\pi \in \LP(n)}
\big( \Psi_n(\pi) - \Psi_n(R \pi) \big) \ket{\pi}
=
{\bf 0}
\ef,
\ee
or also
\be
\label{eq.654898768a1}
\sum_{\pi \in \LP(n)}
\Psi_n(\pi)
\big(
\ket{\pi}
-
\ket{R^{-1} \pi}
\big) 
=
{\bf 0}
\ef.
\ee
If it is understood that a certain operator $\hat{X}$ acts on $\LP(n)$ as
$\hat{X} \ket{\pi} = \ket{X\pi}$, then we do not need to write sums all
the time.
For statements concerning the refined enumerations of FPL,
we will just define, once and forever, the state
\begin{align}
\label{eq.defkets}
\ket{s_n} 
&:= 
\sum_{\phi \in \fpl(n,+)} \ket{\pi(\phi)}
=
\sum_{\pi \in \LP(n)} \Psi_n(\pi) \ket{\pi}
\ef,
\end{align}
and, for example, the dihedral symmetry reads in these notations
\be
(R-1) \ket{s_n} = {\bf 0}
\ef.
\ee
Note that the Temperley-Lieb operators $e_j$, defined in the previous
section, act on $\LP(n)$, and thus expressions such as $e_j \ket{s_n}$
make sense in this notational framework:
\be
e_j \ket{s_n}
=
e_j
\Big(
\sum_{\pi \in \LP(n)} \Psi_n(\pi) \ket{\pi}
\Big)
=
\sum_{\pi \in \LP(n)} \Psi_n(\pi) \ket{e_j (\pi)}
\ef.
\ee
Similarly,
we may have operators $\widetilde{X}$
acting on $\fpl(n,\tau)$.
For dealing with these cases, we will introduce a vector space whose
basis vectors are all the valid FPL configurations, $\kket{\phi} \in
\fpl(n,\tau)$.\footnote{In order to improve readability,
  we use double parenthesis for vectors in this different space.}
This space is thus isomorphic to $\mathbb{C}^{\fpl(n,\tau)}$, 
and the action on the basis vectors is just
\be
\widetilde{X} \kket{\phi} = \kket{ \widetilde{X}(\phi) }
\ef.
\ee
We have natural maps
$\Pi_{\tau}: \fpl(n,\tau) \to \LP(n)$, defined as
$\Pi_{\tau} \kket{\phi} = \ket{\pi(\phi)}$, and also a natural
definition of the states enumerating all FPL
\be
\label{eq.defkets3}
\kket{s_{n,\tau}} 
:= 
\sum_{\phi \in \fpl(n,\tau)} \kket{\phi}
\ee
such that, in particular, according to our definition
(\ref{eq.defkets}) of the state $\ket{s}$,
\be
\ket{s_n} = \Pi_{\pm} \kket{s_{n,\pm}} 
\ef,
\ee
and we could be interested, e.g., in the action
\be
\Pi_+
\widetilde{X}
\kket{s_{n,+}} 
= 
\Pi_+
\sum_{\phi \in \fpl(n,+)}
\kket{ \widetilde{X}(\phi) }
=
\sum_{\phi \in \fpl(n,+)}
\ket{ \pi(\widetilde{X}(\phi)) }
\ef,
\ee
which is a certain vector in $\mathbb{C}^{\LP(n)}$, and thus, for
example, is comparable to $\ket{s_n}$, our vector of interest.

\subsection{The conjecture}

\noindent
Consider FPL configurations in the ensemble $\fpl(n;+)$, and the
Temperley-Lieb Algebra with $2n$ generators.
Adopt the notation $\ket{s_n}$ as in (\ref{eq.defkets}).
Define the \emph{Hamiltonian} (a term motivated by the XXZ Spin Chain)
\be
\label{eq.defH}
H_n=\sum_{k=1}^{2n} e_k
\ef.
\ee
The Razumov-Stroganov conjecture reads
\begin{conj}[Razumov-Stroganov]
\label{conj.RS}
\be
\label{eq.RS1}
H_n \ket{s_n} = 2n \ket{s_n}
\ef.
\ee
\end{conj}
In order to have simple notations, we define
\be
\label{eq.defRS}
\RS_n := (H_n-2n) \ket{s_n}
\ef,
\ee
and the conjecture just states that $\RS_n={\bf 0}$ as a vector in
the linear space~$\mathbb{C}^{\LP(n)}$.



\section{Proof of the conjecture}
\label{sec.proof}

\subsection{A rewriting of the quantity $H \ket{s}$}

\noindent
In a sequence of sections, we analyse the Razumov-Stroganov conjecture
for the periodic $O(1)$ loop model with $2n$ sites, corresponding to
FPL configurations over a square domain of side $n$.
Subscripts $n$, such as in equations (\ref{eq.defkets}) and
(\ref{eq.defH}-\ref{eq.defRS}), will be dropped from now on, in order
to enlight notation.

Choose to fix the boundary conditions, and the labels of the external
black legs, in such a way that the vertical external edge at the
bottom-left corner is black and has label 1, and the labels are given
cyclically in counter-clockwise order (cfr.~figure \ref{fig.stateS} in
\ref{app.gloss}).

The property (\ref{eq.34657865}) allows to rewrite the Hamiltonian
(\ref{eq.defH}) as
\begin{align}
\label{eq.4327658}
H 
&= 
\sum_{k=0}^{2n-1}
R^{k} e_j R^{-k}
\ef,
\end{align}
for any index $1 \leq j \leq 2n$.
Recall that the ordinary Wieland Theorem on dihedral symmetry gives
\be
\label{eq.dih}
R \ket{s} = \ket{s}
\ee
which, combined with (\ref{eq.4327658}), gives
\be
H \ket{s}
=
(1+R+R^{2}+\cdots+R^{2n-1}) 
e_j
\ket{s}
\ef.
\ee
Call $\Sym$ the operator
\be
\Sym = \sum_{k=0}^{2n-1} R^k
\ef,
\ee
which has the simple property
\be
\Sym \; R = R \; \Sym = \Sym
\ef.
\ee
This gives a rewriting of the quantity appearing in the conjecture
\be
\label{eq.H4}
\RS
=
\Sym \,
\big( e_j - 1 \big)
\ket{s}
\ef,
\ee
for any $1 \leq j \leq 2n$.


\subsection{A restatement of the conjecture}

\noindent
For a given plaquette $\alpha$ within our domain $\Lambda$, define the
operator $\nn_{\alpha}(\phi)$ as
\be
\label{eq.nnalph}
\nn_{\alpha}(\phi) =
\left\{
\begin{array}{ll}
+1 & \phi|_{\alpha} = \sqqH \\
-1 & \phi|_{\alpha} = \sqqV \\
 0 & \textrm{otherwise}
\end{array}
\right.
\ee
By ``$\phi|_{\alpha} = \sqqH$'' we mean that the plaquette $\alpha$ is
composed of two black horizontal edges, and two white vertical edges,
while by ``$\phi|_{\alpha} = \sqqV$'' we mean the analogous statement with
black and white interchanged.
Similarly define the operator $\widetilde{\nn}_{\alpha}$, acting
diagonally over
$\mathbb{C}^{\fpl(\Lambda, \tau)}$ as
$\widetilde{\nn}_{\alpha} \kket{\phi} = \nn_{\alpha}(\phi)
\, \kket{\phi}$.

As $\nn_{\alpha}(\phi)\in \{-1,0,+1\}$, the operator 
$\widetilde{\nn}_{\alpha}$ is the difference of two orthogonal
projectors, more precisely, for a state 
$\kket{s} \in \mathbb{C}^{\fpl(\Lambda,\tau)}$,
the state $\widetilde{\nn}_{\alpha} \kket{s}$ takes the form
\begin{align}
\label{eq.NtilGen}
\kket{s}
&=
\hspace{-2.1mm}
\raisebox{-25pt}{
\setlength{\unitlength}{10pt}
\begin{picture}(6,6)
\put(0,0){\includegraphics[scale=2]{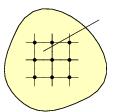}}
\put(6,5.3){$\alpha$}
\end{picture}
}
&
\widetilde{\nn}_{\alpha} \kket{s}
&=
\raisebox{-25pt}{\includegraphics[scale=2]{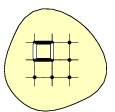}}
-
\raisebox{-25pt}{\includegraphics[scale=2]{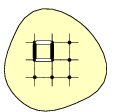}}
\end{align}
In \cite[sec.~5]{wie} it is explained that, for the square domain, two
gyration operations, $H_{\pm}$, can be defined, and the full gyration
operator, by which one proves the symmetry statement
(\ref{eq.763765876}), is $G=H_- H_+$.  In Section \ref{sec.wie}, we
illustrate under which conditions this fact extends to other domains
$(\Lambda, \tau)$.

Assume here that, for such a domain, two gyration
operations $H_{\pm}: \fpl(\Lambda, \tau) \leftrightarrow \fpl(\Lambda,
\bar{\tau})$ are defined, and call $G=H_- H_+$.
For a configuration $\phi \in \fpl(\Lambda, \tau)$, call
$\fpl(\Lambda,\tau;\mathcal{O}(\phi))$ 
the orbit of $\phi$ inside $\fpl(\Lambda, \tau)$,
under the action of $G$. 

We have the lemma
\begin{lemma}
\label{thm.HVgen}
With the definitions above, for every plaquette $\alpha \in \Lambda$,
and every $\phi \in \fpl(\Lambda, \tau)$,
we have
\be
\sum_{\phi' \in \fpl(\Lambda,\tau;\mathcal{O}(\phi)) }
\nn_{\alpha}(\phi')
= 0
\ef.
\ee
\end{lemma}
Now consider 
the $n \times n$ square domain with
alternating boundary conditions, ensemble $\fpl(n,\pm)$. 

\begin{figure}[tb!]
\begin{center}
\includegraphics[scale=1.25]{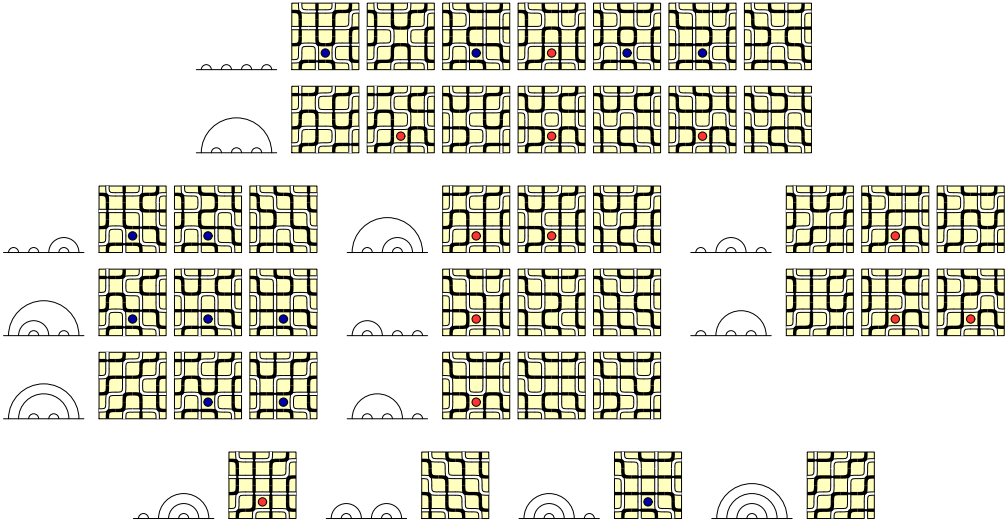}
\end{center}
\caption{An illustration of Corollary \ref{cor.HV}, on the 42 FPL in
  the square of side 4. The three blocks of the picture (first two
  rows, rows 3--5 and the last row) correspond to the three classes
  $[\pi]$ in $\LP^*(4)$. We have blue and red bullets for
  configurations corresponding to the two non-trivial cases of
  equation (\ref{eq.nnalph}), for the plaquette at coordinate $(3,2)$.
  As claimed, separately in each of the three blocks, there are as
  many blue bullets as red bullets.\label{fig.cor.HV}}
\end{figure}

Call $\LP^*(n)$ the set of link patterns in $\LP(n)$, quotiented
w.r.t.\ cyclic rotations, and call $[\pi]$ an element in this set (the
class of $\pi$ w.r.t.\ the equivalence relation $\pi \sim \pi'$ iff
$\pi=R^k \pi'$ for some $0 \leq k \leq 2n-1$).  Call
$\fpl(n,\pm;[\pi])$ the refined subsets of $\fpl(n,\pm)$ w.r.t.\ the
quantities $[ \pi(\phi) ]$.  Because of the dihedral symmetry for FPL
on the square domain, the sets $\fpl(n,\pm;[\pi])$ are a disjoint
union of whole orbits $\fpl(\Lambda,\tau;\mathcal{O}(\phi))$, so we
get the corollary
\begin{corollary}
\label{cor.HV}
For any $[\pi] \in \LP^*(n)$, and any plaquette $\alpha$
\begin{align}
\label{eq.incor.HV}
\sum_{\phi \in \fpl(n,\pm;[\pi])}
\nn_{\alpha}(\phi)
= 0
\ef.
\end{align}
\end{corollary}
In vector notation, (\ref{eq.incor.HV}) is equivalent to
\begin{align}
\label{eq.24687634a}
\Sym\;
\Pi_{\pm}
\widetilde{\nn}_{\alpha}
\kket{s_{n,\pm}}
= 0
\ef.
\end{align}
Indeed,
\be
\begin{split}
\Sym\;
\Pi_{\pm}
\widetilde{\nn}_{\alpha}
\kket{s_{n,\pm}}
&=
\Sym\;
\Pi_{\pm}
\sum_{\phi \in \fpl(n,\pm)}
\widetilde{\nn}_{\alpha}
\kket{\phi}
=
\Sym\;
\Pi_{\pm}
\sum_{\phi \in \fpl(n,\pm)}
\nn_{\alpha}(\phi)
\kket{\phi}
\\
&=
\Sym\;
\sum_{\phi \in \fpl(n,\pm)}
\nn_{\alpha}(\phi)
\ket{\pi(\phi)}
\\
&=
\sum_{\pi \in \LP(n)}
| \mathrm{Aut}(\pi) |
\Big(
\sum_{\phi \in \fpl(n,\pm;[\pi])}
\nn_{\alpha}(\phi)
\Big)
\ket{\pi}
\ef,
\end{split}
\ee
where $| \mathrm{Aut}(\pi) |$ is the cardinality of the subgroup of
rotations that stabilize~$\pi$.
This corollary is pictorially illustrated in an example in
figure~\ref{fig.cor.HV}.

Lemma \ref{thm.HVgen} and Corollary~\ref{cor.HV} are interesting by
themselves. However, at this point we prefer to stress immediately
what will show up to be their crucial property
\begin{proposition}
\label{prop.equiv}
For the $n \times n$ square, call $\alpha_j$ the plaquette located in
the $(2j-1)$-th column-position along the bottom row.  The quantity in
the Razumov-Stroganov conjecture, defined in (\ref{eq.defRS}), is
equal to
\be
\label{eq.24687634}
\RS_n = 
\sum_{j=1}^{\lceil n/2 \rceil} 
\Sym\;
\Pi_{\pm}
\widetilde{\nn}_{\alpha_j}
\kket{s_{n,\pm}}
\ef.
\ee
\end{proposition}
Clearly, the Razumov-Stroganov conjecture, equation
(\ref{eq.RS1}),
is proven if both Lemma \ref{thm.HVgen} and Proposition \ref{prop.equiv} are
proven, as the right-hand side of (\ref{eq.24687634}) is a sum of
quantities as in (\ref{eq.24687634a}), that vanish as a result of
Corollary~\ref{cor.HV}.

We give here the proof of Lemma \ref{thm.HVgen}, which is relatively
short and simple, and devote the rest of the paper to the more
composite proof of Proposition~\ref{prop.equiv}.

The proof of Lemma \ref{thm.HVgen} can be read at two levels. The
reader aware of \cite{wie}, and interested only in the case of the
square, $\fpl(\Lambda,\tau) = \fpl(n,+)$, sufficient at the purpose of
Corollary~\ref{cor.HV} and thus of the Razumov-Stroganov conjecture,
can read directly the proof, with the understanding of this
restriction. The reader interested in the more general statement will
find in Section \ref{ssec.wie1} the required preliminary discussion on
gyration.

\vspace{2mm}
\noindent
{\it Proof of Lemma \ref{thm.HVgen}.}  Call $\tilde{H}_+$ and
$\tilde{H}_-$ respectively the maps $H_+$ and $H_-$, followed by
complementation. As the complementation commutes with $H_{\pm}$ (it is
just a labeling of the colours, and the definition (\ref{eq.defHphi})
is symmetric), we also have $G= \tilde{H}_- \tilde{H}_+$, and all
these three maps are bijections over $\fpl(\Lambda, \tau)$ (without
involving $\bar{\tau}$).

We start by proving the statement for plaquettes $\alpha$ adjacent to
the border.
For any $\phi \in \fpl(n,\pm)$, define the infinite string in the alphabet
$\{-1,0,+1\}$
\be
\bm{\nu}(\phi) = \big( \nn_{\alpha}(\phi), \nn_{\alpha}(G \phi), 
\nn_{\alpha}(G^2 \phi), \nn_{\alpha}(G^3 \phi), \ldots \big)
\ef.
\ee
As the string of $\{ G^k \phi \}_{k \in \mathbb{N}}$ is periodic over the orbit
$\mathcal{O}(\phi)$, also the values $\nn_{\alpha}(G^k \phi)$ are a
periodic sequence, with period $|\mathcal{O}(\phi)|$.
We claim that $\bm{\nu}(\phi)$ is composed of alternating $+1$'s and
$-1$'s, separated by intervals of zeroes (possibly empty). From this
statement, 
the proposition specialized to plaquettes on the border would follow.

The analysis of $\bm{\nu}(\phi)$ is performed through the analysis of
a further auxiliary string. Assume that $\alpha$ is adjacent to the
border of the square through a horizontal edge $e$, and that it
undergoes gyration in the first of the two parity rounds,
$\tilde{H}_+$.\footnote{In the other cases the reasoning would be
  modified in a minor way. However, curiously, within
  Proposition~\ref{prop.equiv} we only need this case, and only in the
  easier case of plaquettes on the border.}  
Define the string in the alphabet $\{b,w\}$
\be
\bm{\mu}(\phi) = \big( \phi_e, (G \phi)_e, 
(G^2 \phi)_e, (G^3 \phi)_e, \ldots \big)
\ef.
\ee
Then, clearly $\nu_k=+1$ only if $\mu_k=b$, and $\nu_k=-1$ only if
$\mu_k=w$.  But we have more than this.  Indeed, if $\nu_k=+1$, the
plaquette $\alpha$ will undergo gyration in the next $\tilde{H}_+$
round. Then, $e$ is not touched by the $\tilde{H}_-$ round. So we have
not only $\mu_k=b$, but also $\mu_{k+1}=w$. Similarly, if $\nu_k=-1$,
we have not only $\mu_k=w$, but also $\mu_{k+1}=b$. Also the converse holds:
the only possibility for $\mu_k \neq \mu_{k+1}$ is that $\nu_k \neq
0$, as otherwise $e$ is not interested by gyration at round
$\tilde{H}_+$, and in general $e$ is never interested by
$\tilde{H}_-$. So, the sequence $\bm{\nu}(\phi)$ collects the
positions of the inversions (with sign) in the binary sequence
$\bm{\mu}(\phi)$, and thus has the claimed structure of an alternating
sequence of $+1$'s and $-1$'s, separated by intervals of zeroes. This
completes the proof for boundary plaquettes.

Now consider two neighbouring plaquettes $\alpha$ and $\beta$, sharing
a common edge $e$, say horizontal. We claim that
\be
\label{eq.67647543754}
\sum_{\phi' \in \fpl(\Lambda,\tau;\mathcal{O}(\phi)) }
\nn_{\alpha}(\phi')
\ -
\sum_{\phi' \in \fpl(\Lambda,\tau;\mathcal{O}(\tilde{H}_+(\phi))) }
\nn_{\beta}(\phi')
= 0
\ef.
\ee
From this statement, the whole lemma would follow, as we already
know that, if $\alpha$ is a border plaquette, in the equation above
the sum on the left is zero.


The reasoning is analogous to the previous one, but now we consider a
string $\tilde{\bm{\mu}}$ on the single rounds $\tilde{H}_{\pm}$
\be
\tilde{\bm{\mu}}(\phi) = \big( \phi_e, (\tilde{H}_+ \phi)_e, 
(G \phi)_e, (\tilde{H}_+ G \phi)_e, 
(G^2 \phi)_e, (\tilde{H}_+ G^2 \phi)_e, 
\ldots \big)
\ef.
\ee
We also consider the strings
\begin{align}
\bm{\nu}_{\alpha}(\phi) 
&= \big( \nn_{\alpha}(\phi), \nn_{\alpha}(G \phi), 
\nn_{\alpha}(G^2 \phi), \nn_{\alpha}(G^3 \phi), \ldots \big)
\ef;
\\
\bm{\nu}_{\beta}(\phi) 
&= \big( \nn_{\beta}(\tilde{H}_+ \phi), \nn_{\beta}(\tilde{H}_+ G \phi), 
\nn_{\beta}(\tilde{H}_+ G^2 \phi), \nn_{\beta}(\tilde{H}_+ G^3 \phi), \ldots \big)
\ef.
\end{align}
A typical example could be
\[
\begin{array}{r|ccccccccccccccl}
\hline
\bm{\nu}_{\alpha} & 0 &   & + &   & - &   & - &   & 0 &   & - &   & + &   & \cdots \\
\tilde{\bm{\mu}}  & b & b & b & w & w & b & w & b & w & w & w & b & b & w & \cdots \\
\bm{\nu}_{\beta}  &   & 0 &   & 0 &   & + &   & + &   & 0 & & 0 &   & - & \cdots \\
\hline
\end{array}
\]
An argument completely analogous to the one exploited in the
border-plaquette case shows that the inversions in the string
$\tilde{\bm{\mu}}$ ($b \to w$ and $w \to b$) are in correspondence
with the positions of $+1$ and $-1$ along the strings
$\bm{\nu}_{\alpha}$ and $\bm{\nu}_{\beta}$ ($+1$ for $b \to w$, $-1$
for $w \to b$, and along the $\alpha$ or $\beta$ string depending on the
parity of the position of the inversion along~$\tilde{\bm{\mu}}$).

While the string $\bm{\nu}_{\alpha}$ is just circuitating along the
orbit of $\phi$ w.r.t.\ the action of $G$, the string
$\bm{\nu}_{\beta}$ is circuitating along the orbit of $H_+(\phi)$
w.r.t.\ the action of $\tilde{H}_+ \tilde{H}_-$, which is \emph{not}
$G$. However, it is $G^{-1}$, up to a conjugation with complementation
(i.e., it is $\overline{G^{-1}(\bar{\phi})}$). But $G$ and $G^{-1}$
have the same orbit, and the complementation relates the orbit
$\mathcal{O}(\phi))$ over $\fpl(\Lambda,\tau)$ to the orbit
$\mathcal{O}(\bar{\phi}))$ over $\fpl(\Lambda,\bar{\tau})$, which have
opposite sets of $\nn_{\beta}(\phi')$ values.  This proves our claim
(\ref{eq.67647543754}) (justifying the minus sign), and completes the
proof of the lemma.  \hfill $\square$

\subsection{Definition of auxiliary combinations}
\label{ssec.auxil}

\noindent
Call $\nb = \lceil n/2 \rceil$ and $\nc =\lfloor n/2 \rfloor$.  For $1
\leq j \leq \nb$, call $\ssb{a,b,c}{j}$ the state over the $n \times
n$ square, with enumerations of FPL in the ensemble $\fpl(n;+)$,
restricted to the case in which the $(2j-1)$-th node of the last row
is an $a$, $b$ or $c$ configuration
(w.r.t.\ the definition (\ref{eq.defabc})).  For $1 \leq j \leq \nc$, call
$\ssc{a,b,c}{j}$ the state with enumerations of FPL in the ensemble
$\fpl(L;+)$, restricted to the case in which the $(2j)$-th node of the
last row is an $a$, $b$ or $c$ configuration.

These combinations, and various others that we will need along the
proof, are illustrated in a glossary in \ref{app.gloss}.

The resulting domains, restricted by a single site, have in general
some frozen regions, i.e.\ regions of the square domain in which the
configuration is fixed in any valid FPL, and we can read the states
above as states over smaller domains. 

States $\ssb{c}{j}$ and $\ssc{c}{j}$ force restriction over the whole
last row, and furthermore, for $n$ even, $\ssb{c}{1}$ and
$\ssc{c}{\nc}$ also force restriction respectively over the first and
last column, leading to FPL configurations over the smaller $(n-1)
\times (n-1)$ square domain. Similarly, for $n$ odd, we have this
property for $\ssb{c}{1}$ and $\ssb{c}{\nb}$.  However, in what
follows we shall \emph{not} need these last properties (and our proof
is not inductive).

States $\ssb{b}{j}$ and $\ssc{b}{j}$ force restriction over the part
of the
last row which is on the left of the decimated site, while
states $\ssb{a}{j}$ and $\ssc{a}{j}$ force restriction over the part of
the last row which is on the right. We have, in particular,
\begin{subequations}
\label{eqs.sabzero}
\begin{align}
\ssb{a}{1} &= \ssc{b}{\nc} = 0
&&
\textrm{$n$ even;}
\\
\ssb{a}{1} &= \ssb{b}{\nb} = 0
&&
\textrm{$n$ odd;}
\end{align}
\end{subequations}
as these choices of restriction on the corner sites are inconsistent with
the boundary conditions.

\subsection{Identities}

\noindent
Any valid FPL configuration has exactly one $c$ entry in the last
row. More precisely, the entries in order in the last row have the
form
\[
(b,b,\ldots,b,
\!\!\stackrel{\raisebox{3pt}{\scriptsize{$i$-th}}}{c}\!\!,
a,a,\ldots,a)
\]
for some $1 \leq i \leq n$.  This leads to a refinement of the
enumerations
\begin{proposition}[Last-row decomposition]
\label{prop.LRD}
\be
\ket{s}
=
\sum_{j=1}^{\nb}
\ssb{c}{j}
+
\sum_{j=1}^{\nc}
\ssc{c}{j}
\ef.
\ee
\end{proposition}
Furthermore, we can refine the enumerations w.r.t.\ the three choices
among $a$, $b$, $c$ for any single site in the last row, getting
\begin{proposition}[One-site expansion]
\label{prop.OSE}
\begin{align}
\label{eq.onesiteB}
\ket{s}
&=
\ssb{a}{j}
+
\ssb{b}{j}
+
\ssb{c}{j}
&&
\forall \quad 1 \leq j \leq \nb
\ef;
\\
\label{eq.onesiteC}
\ket{s}
&=
\ssc{a}{j}
+
\ssc{b}{j}
+
\ssc{c}{j}
&&
\forall \quad 1 \leq j \leq \nc
\ef.
\end{align}
\end{proposition}
We have the simple fact 
\be
\label{eq.eceqc}
e_j \ssc{c}{j} = \ssc{c}{j}
\ee
as the corresponding restriction forces an arc between $j$ and $j+1$,
already within the last row, which is frozen
(cfr.\ figure~\ref{fig.stateSc}, right column).

We have simple \emph{recursion relations} for $\ssc{a}{j}$,
$\ssc{b}{j}$, $\ssb{a}{j}$ and $\ssb{b}{j}$ states, performed by further
refining the configurations over another site. For example, the state
$\ssc{a}{j}$ is already restricted to the $(2j)$-th site of the last
row to be an $a$, thus, from the $(2j+1)$-th site on, the row is frozen
to be filled with $a$'s, and on the $(2j-1)$-th site we can find
only either an $a$ or a $c$, in the two cases corresponding to the
classes $\ssb{a}{j}$ and $\ssb{c}{j}$ respectively.
Reasonings in this fashion lead to the set of equations
\begin{subequations}
\begin{align}
\label{eq.ricb2}
\ssb{a}{j} 
&= 
\ssc{c}{j-1}
+
\ssc{a}{j-1}
\ef;
\\
\label{eq.rica1}
\ssb{b}{j} 
&= 
\ssc{c}{j}
+
\ssc{b}{j}
\ef;
\\
\label{eq.ricb1}
\ssc{a}{j} 
&= 
\ssb{c}{j}
+
\ssb{a}{j}
\ef;
\\
\label{eq.rica2}
\ssc{b}{j} 
&= 
\ssb{c}{j+1}
+
\ssb{b}{j+1}
\ef.
\end{align}
\end{subequations}
In Section \ref{sec.wie}, we generalize the analysis of the gyration
operation, performed by Wieland in \cite{wie}, to arbitrary regions of
the square lattice, and arbitrary boundary conditions. This analysis,
specialized to our states $\ket{s^{a,b,c}_j}$, leads to a number of
relations that will have a crucial role in what follows.
\begin{proposition}[Gyration relations]
\label{prop.gyrat}
\begin{align}
\label{eq.2354765aa}
e_{j} \ssb{a}{j} 
&= 
R^{-1} e_{j-1} \ssb{a}{j}
\ef;
\\
e_j \ssb{b}{j} 
&= 
R^{-1} e_{j-1} \ssb{b}{j}
\ef;
\\
\label{eq.2354765a}
e_j \ssb{c}{j} 
&= 
R^{-1} e_{j-1} \ssb{c}{j}
\ef;
\\
\label{eq.2354765aaa}
e_j \ssc{a}{j} 
&= 
R^{-1} e_{j-1} \ssc{a}{j}
\ef;
\\
e_j \ssc{b}{j} 
&= 
R e_{j+1} \ssc{b}{j}
\ef.
\end{align}
\end{proposition}
(The proof of these relations is postponed to Section
\ref{sec.gyraproofs}).
We do not have a relation for $\ssc{c}{j}$, as we have instead the
stronger and easier fact~(\ref{eq.eceqc}).

A further relation holds in very general circumstances. Take a generic
domain $\Lambda$ for FPL configurations, and consider the two boundary
conditions $\tau_1 = (x_1,\ldots,x_{2n-2},b,b)$, 
$\tau_2 = (x_1,\ldots,x_{2n-2},w,w)$. Also assume that the two sites
$v$, $v'$ adjacent to the two last external legs are connected in
$\Lambda$ by an edge $e$. Say that $\tau_1$ has $2m$ black legs (so
$\tau_2$ has $2m-2$ black edges).

\begin{figure}
\begin{align*}
&\rule{0pt}{45pt}
\setlength{\unitlength}{10pt}
\begin{picture}(6,6.5)(0,-.5)
\put(0.5,5.5){\makebox[0pt][r]{$\ket{s_1}$}}
\put(0,0){\includegraphics[scale=2]{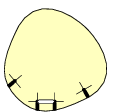}}
\put(-.05,-.9){$\scriptstyle{2m-1}$}
\put(2.9,-.9){$\scriptstyle{2m}$}
\put(5.2,0){$\scriptstyle{1}$}
\put(6.5,1){$\vdots$}
\put(-2.1,1){$\scriptstyle{2m-2}$}
\end{picture}
&&
\setlength{\unitlength}{10pt}
\begin{picture}(6,6.5)(0,-.5)
\put(0.5,5.5){\makebox[0pt][r]{$\ket{s_2}$}}
\put(0,0){\includegraphics[scale=2]{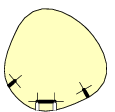}}
\put(5.2,0){$\scriptstyle{1}$}
\put(6.5,1){$\vdots$}
\put(-2.1,1){$\scriptstyle{2m-2}$}
\end{picture}
&&
\setlength{\unitlength}{10pt}
\begin{picture}(6,6.5)(0,-.5)
\put(0.5,5.5){\makebox[0pt][r]{$a_{2m-1} \ket{s_2}$}}
\put(0,0){\includegraphics[scale=2]{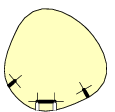}}
\put(-.7,-.9){$\scriptstyle{2m-1}$}
\put(3.8,-.9){$\scriptstyle{2m}$}
\put(5.2,0){$\scriptstyle{1}$}
\put(6.5,1){$\vdots$}
\put(-2.1,1){$\scriptstyle{2m-2}$}
\end{picture}
\end{align*}

\caption{A graphical representation of the three state
$\ket{s_1}$ (left), $\ket{s_2}$ (center) and $a_{2m-1} \ket{s_2}$ (right).
\label{fig.SPR}}
\end{figure}

For $\pi \in \LP(m)$, call $\Psi_1(\pi)$ the number of
configurations $\phi \in \fpl(\Lambda, \tau_1)$, with link pattern
$\pi$, and such that $\phi(e)=w$. For $\pi' \in \LP(m-1)$, call
$\Psi_2(\pi')$ the number of configurations $\phi \in \fpl(\Lambda,
\tau_2)$, with link pattern $\pi'$, and such that $\phi(e)=b$. 
Call $a_{2m-1}$ a map from $\LP(m-1)$ to $\LP(m)$
which adjoin the arc $(2m-1,2m)$ to the link pattern (we will come
back to this sort of operators in Section \ref{sec.TLca}).
Define the states in $\mathbb{C}^{\LP(m)}$
and $\mathbb{C}^{\LP(m-1)}$
\begin{align}
\ket{s_1}
&= \sum_{\pi \in \LP(m)} \Psi_1(\pi) \,\ket{\pi}
\ef;
&
\ket{s_2}
&= \sum_{\pi' \in \LP(m-1)} \Psi_2(\pi') \,\ket{\pi'}
\ef.
\end{align}
A graphical illustration of the three states
$\ket{s_1}$, $\ket{s_1}$ and $a_{2m-1} \ket{s_1}$
is in figure~\ref{fig.SPR}

Then we have
\begin{proposition}[Simple path reversal]
\label{prop.SPR}
\be
e_{2 m-1} \ket{s_1} = a_{2m-1} \ket{s_2}
\ef.
\ee
\end{proposition}
{\it Proof.}  First remark that the configurations $\phi$ contributing
to the two ensembles are in bijection, corresponding to reverse the
values of $\phi$ over the three-edge path including $e$ and the two
boundary edges adjacent to $v$ and $v'$.

Now consider a configuration $\phi$ contributing to $\Psi_1(\pi)$. If
we drop the two boundary black edges adjacent to $v$ and $v'$, these
sites are endpoints of black open paths, and thus are either connected
one with the other, or to some two black endpoints, with indices $j$,
$j'$.  In the associated configuration $\phi'$, we get in the two
cases respectively that $v$ and $v'$ are part of a cycle, and that we
have an arc connecting $j$ and $j'$, passing through $v$ and $v'$. All
the rest of the configuration is unperturbed, and thus also all the
rest of the link pattern, so, in both cases $\phi'$ contributes to
$\Psi_2(\pi')$, for the only $\pi' \in \LP(m-1)$ such that
$a_{2m-1} (\pi') = e_{2m-1} (\pi)$.
\hfill $\square$

\subsection{Proof of the equivalence statement}

\noindent
Now we have all the ingredients for proving Proposition
\ref{prop.equiv}, and we can start studying the quantity $\RS$ in its
form of (\ref{eq.H4}).  We thus analyse $e_j \ket{s}$, for a whatever
$1 \leq j \leq \nb$.\footnote{We choose this range because, for
  simplicity, we only introduced auxiliary states $\ket{s^{a,b,c}_j}$
  for refinements along the bottom line. The proof could be done with
  minor modifications for any $1 \leq j \leq 2n$.}
Do the $j$-th one-site expansion (\ref{eq.onesiteB})
(this coincidence of indices is a key point)
\be
\begin{split}
e_j \ket{s}
&=
e_j 
\big(
\ssb{a}{j}
+
\ssb{b}{j}
+
\ssb{c}{j}
\big)
\\
&=
\ssc{c}{j}
+
e_j 
\ssb{a}{j}
+
e_j 
\ssc{b}{j}
+
e_j 
\ssb{c}{j}
\ef.
\end{split}
\ee
where we used (\ref{eq.eceqc}) and (\ref{eq.rica1}).  We could have
similarly done the $j$-th one-site expansion (\ref{eq.onesiteC}), and
obtain the same result, by using (\ref{eq.eceqc}) and
(\ref{eq.ricb1}).

\begin{figure}[tb]
\[
\begin{array}{cc}
\ssbcx{j}
&
e_{j} \ssbcx{j}
\\
\setlength{\unitlength}{10pt}
\begin{picture}(14,4.5)(-.5,0)
\put(0.25,0.5){\includegraphics[scale=2]{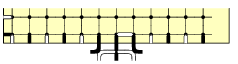}}
\put(-0.2,2.1){\makebox[0pt][c]{$\scriptstyle{2n}$}}
\put(2,1){\makebox[0pt][c]{$\scriptstyle{1}$}}
\put(4,1){\makebox[0pt][c]{$\scriptstyle{2}$}}
\put(5,1){\makebox[0pt][c]{$\scriptstyle{\cdots}$}}
\put(5.2,0.3){\makebox[0pt][c]{$\scriptstyle{j-1}$}}
\put(7,0){\makebox[0pt][c]{$\scriptstyle{j}$}}
\put(8.8,0.3){\makebox[0pt][c]{$\scriptstyle{j+1}$}}
\put(10,1){\makebox[0pt][c]{$\scriptstyle{j+2}$}}
\put(11,1){\makebox[0pt][c]{$\scriptstyle{\cdots}$}}
\end{picture}
&
\setlength{\unitlength}{10pt}
\begin{picture}(14,4.5)(-.5,0)
\put(0.25,0.5){\includegraphics[scale=2]{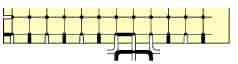}}
\put(-0.2,2.1){\makebox[0pt][c]{$\scriptstyle{2n}$}}
\put(2,1){\makebox[0pt][c]{$\scriptstyle{1}$}}
\put(4,1){\makebox[0pt][c]{$\scriptstyle{2}$}}
\put(4.7,1){\makebox[0pt][c]{$\scriptstyle{\cdots}$}}
\put(5.7,1){\makebox[0pt][c]{$\scriptstyle{j-1}$}}
\put(7,0){\makebox[0pt][c]{$\scriptstyle{j}$}}
\put(9,0){\makebox[0pt][c]{$\scriptstyle{j+1}$}}
\put(10.3,1){\makebox[0pt][c]{$\scriptstyle{j+2}$}}
\put(11.5,1){\makebox[0pt][c]{$\scriptstyle{\cdots}$}}
\end{picture}
\\
\setlength{\unitlength}{10pt}
\begin{picture}(14,4.5)(-.5,0)
\put(0.25,0.5){\includegraphics[scale=2]{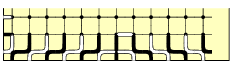}}
\put(-0.2,2.1){\makebox[0pt][c]{$\scriptstyle{2n}$}}
\put(1,0){\makebox[0pt][c]{$\scriptstyle{1}$}}
\put(3,0){\makebox[0pt][c]{$\scriptstyle{2}$}}
\put(3.8,0){\makebox[0pt][c]{$\scriptstyle{\cdots}$}}
\put(5,0){\makebox[0pt][c]{$\scriptstyle{j-1}$}}
\put(7,0){\makebox[0pt][c]{$\scriptstyle{j}$}}
\put(9,0){\makebox[0pt][c]{$\scriptstyle{j+1}$}}
\end{picture}
&
\setlength{\unitlength}{10pt}
\begin{picture}(14,4.5)(-.5,0)
\put(0.25,0.5){\includegraphics[scale=2]{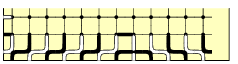}}
\put(-0.2,2.1){\makebox[0pt][c]{$\scriptstyle{2n}$}}
\put(1,0){\makebox[0pt][c]{$\scriptstyle{1}$}}
\put(3,0){\makebox[0pt][c]{$\scriptstyle{2}$}}
\put(3.8,0){\makebox[0pt][c]{$\scriptstyle{\cdots}$}}
\put(5,0){\makebox[0pt][c]{$\scriptstyle{j-1}$}}
\put(7,0){\makebox[0pt][c]{$\scriptstyle{j}$}}
\put(9,0){\makebox[0pt][c]{$\scriptstyle{j+1}$}}
\put(11,0){\makebox[0pt][c]{$\scriptstyle{j+2}$}}
\end{picture}
\\
\setlength{\unitlength}{10pt}
\begin{picture}(14,4.5)(-.5,0)
\put(0.25,0.5){\includegraphics[scale=2]{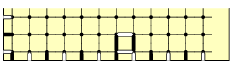}}
\put(-0.2,2.1){\makebox[0pt][c]{$\scriptstyle{2n}$}}
\put(1,0){\makebox[0pt][c]{$\scriptstyle{1}$}}
\put(3,0){\makebox[0pt][c]{$\scriptstyle{2}$}}
\put(3.8,0){\makebox[0pt][c]{$\scriptstyle{\cdots}$}}
\put(5,0){\makebox[0pt][c]{$\scriptstyle{j-1}$}}
\put(7,0){\makebox[0pt][c]{$\scriptstyle{j}$}}
\put(9,0){\makebox[0pt][c]{$\scriptstyle{j+1}$}}
\end{picture}
&
\setlength{\unitlength}{10pt}
\begin{picture}(14,4.5)(-.5,0)
\put(0.25,0.5){\includegraphics[scale=2]{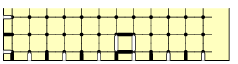}}
\put(-0.2,2.1){\makebox[0pt][c]{$\scriptstyle{2n}$}}
\put(1,0){\makebox[0pt][c]{$\scriptstyle{1}$}}
\put(3,0){\makebox[0pt][c]{$\scriptstyle{2}$}}
\put(3.8,0){\makebox[0pt][c]{$\scriptstyle{\cdots}$}}
\put(5,0){\makebox[0pt][c]{$\scriptstyle{j-1}$}}
\put(7,0){\makebox[0pt][c]{$\scriptstyle{j}$}}
\put(9,0){\makebox[0pt][c]{$\scriptstyle{j+1}$}}
\put(11,0){\makebox[0pt][c]{$\scriptstyle{j+2}$}}
\end{picture}
\end{array}
\]
\caption{The states $\ssbcx{j}$ and $e_{j} \ssbcx{j}$, here
  represented three times: top, as they come out from the definition
  $\ssb{c}{j} = \ssb{cb}{j} + \ssbcx{j}$, and the application of the
  simple path reversal relation; middle, reintegrating the frozen
  region in the full $n \times n$ square; bottom, representing a
  minimal set of constrained edges, in the full $n \times n$ square,
  leading to the same frozen region, and making evident the connection
  with the projector $\tilde{\nn}_{\alpha}$ as described in equation
  (\ref{eq.NtilGen}).\label{fig.stateXj}}
\end{figure}

We can use the recursions (\ref{eq.rica2}) and (\ref{eq.ricb2}) in
order to push the $\ssb{a}{j}$ and $\ssc{b}{j}$ states towards their
values at the left and right corners, respectively, that are zero by
(\ref{eqs.sabzero}). In order to preserve the coincidence of indices
in combinations $e_j \ket{s^{a,b,c}_j}$, we will use the gyration
relations of Proposition \ref{prop.gyrat}. Indeed we have
\begin{align}
\begin{split}
e_j \ssb{a}{j}
&=
R^{-1} e_{j-1} \ssb{a}{j}
=
R^{-1} e_{j-1}
(\ssc{c}{j-1}+\ssb{c}{j-1}+\ssb{a}{j-1})
\\&=
R^{-1}
(\ssc{c}{j-1}+ e_{j-1} \ssb{c}{j-1}+ e_{j-1} \ssb{a}{j-1})
\ef;
\end{split}
\\
\begin{split}
e_j \ssc{b}{j}
&=
R e_{j+1} \ssc{b}{j}
=
R e_{j+1} 
(\ssc{c}{j+1}+\ssb{c}{j+1}+\ssc{b}{j+1})
\\&=
R 
(\ssc{c}{j+1}+ e_{j+1} \ssb{c}{j+1}+ e_{j+1} \ssc{b}{j+1})
\ef.
\end{split}
\end{align}
This procedure drops down $\ssc{c}{\ell}$ and $\ssb{c}{\ell}$ states at every
step, that we can easily collect. The result is
\be
e_j \ket{s}
=
\sum_{\ell=1}^{\nb}
R^{\ell-j} e_{\ell} \ssb{c}{\ell}
+
\sum_{\ell=1}^{\nc}
R^{\ell-j} \ssc{c}{\ell}
\ef;
\ee
and thus
\be
\Sym\; e_j \ket{s}
=
\Sym
\Big(
\sum_{\ell=1}^{\nb}
e_{\ell} \ssb{c}{\ell}
+
\sum_{\ell=1}^{\nc}
\ssc{c}{\ell}
\Big)
\ef;
\ee
independently from the choice of $j$, as expected. We have thus one of
the two summands in (\ref{eq.H4}), the other one being just $\Sym\, \ket{s}$.
For this term, using the last-row decomposition, Proposition
\ref{prop.LRD}, we have
\be
\Sym\, \ket{s}
=
\Sym
\Big(
\sum_{\ell=1}^{\nb}
\ssb{c}{\ell}
+
\sum_{\ell=1}^{\nc}
\ssc{c}{\ell}
\Big)
\ef.
\ee
From this, we get an equivalent formulation of $\RS$
\be
\label{eq.765876a}
\RS
=
\Sym
\sum_{j=1}^{\nb}
\big(e_{j}-1 
\big)
\ssb{c}{j}
\ef.
\ee
We investigate these summands, concentrating on the quantity
\be
\big(e_{j}-1 \big)
\ssb{c}{j}
\ef.
\ee
Consider the state $\ssb{c}{j}$, and the site adjacent to the black
external leg labeled as $j$. This site may be in the configuration
$b$, thus forcing the connectivity among leg $j$ and $j+1$, or in
configuration $a$ or $c$. We call respectively 
$\ssb{cb}{j}$ and $\ssbcx{j}$ these two states. So we get
\be
\label{eq.543865768}
\begin{split}
\big(e_{j}-1 \big)
\ssb{c}{j}
&=
\big(e_{j}-1 \big)
\big(
\ssb{cb}{j} + \ssbcx{j}
\big)
=
\big(e_{j}-1 \big)
\ssbcx{j}
\ef.
\end{split}
\ee
The state $\ssbcx{j}$ is shaped as shown in 
figure \ref{fig.stateXj}, top left.
Remark that this domain is suitable for the application of the simple
path reversal relation, Proposition \ref{prop.SPR}. This allows us to
identify the vector $e_{j} \ssbcx{j}$ in $\mathbb{C}^{\LP(n)}$ with
the vector associated to the state described by the domain in
figure \ref{fig.stateXj}, top right (where we also included a black
arc connecting $j$ and $j+1$, disjoint from the domain, in accordance
with the action of $e_j$).

Restoring the frozen last row in both states $\ssbcx{j}$ and $e_{j}
\ssbcx{j}$ of 
figure \ref{fig.stateXj}, top left and right respectively, leads back
to the $n \times n$
square domain, with alternating boundary condition $\tau_+$, and the
edges around a whole plaquette constrained to certain values
(this is shown, in two steps, in the bottom part of figure
\ref{fig.stateXj}).
The plaquette is the one of the last row, and the column $2j-1$, that
in Proposition \ref{prop.equiv} has been called~$\alpha_j$.

This leads us to recognize
\be
\label{eq.765876c}
e_{j} \ssbcx{j} - \ssbcx{j} 
= \Pi_+ \widetilde{\nn}_{\alpha_j} \kket{s_{+}}
\ef.
\ee
Indeed, the reason why we can replace the combination above over the
state $\ssbcx{j}$, collecting \emph{constrained} FPL configurations,
with the full state $\kket{s_{+}}$, is that the operator
$\widetilde{\nn}_{\alpha_j}$ makes zero on all the configurations
$\phi$ whose restriction to $\alpha_j$ does not coincide neither with
the constraint depicted in 
figure \ref{fig.stateXj}, bottom left
(for which it makes $+1$), nor with
the one depicted in 
figure \ref{fig.stateXj}, bottom right
(for which it makes~$-1$), as explained more
in general at equation~(\ref{eq.NtilGen}).

Collecting (\ref{eq.765876a}), (\ref{eq.543865768}) and
(\ref{eq.765876c}), we conclude that
\be
\label{eq.765876d}
\RS
=
\sum_{j=1}^{\nb}
\Sym\;
\Pi_+ \widetilde{\nn}_{\alpha_j} \kket{s_{+}}
\ef,
\ee
thus completing the proof.
\hfill $\square$

\section{Dihedral symmetry and gyration}
\label{sec.wie}

\subsection{A revisitation of Wieland proof}
\label{ssec.wie1}

\noindent
In \cite{wie}, Wieland proves the dihedral symmetry in the enumeration
of FPL classes with given link pattern $\pi$, in the square $n \times
n$ domain with alternating boundary conditions (a fact previously
conjectured by Bosley and Fidkowski, and unpublished). He proves a
more general fact, for a three-time refined enumeration of FPL,
according to the link pattern $\pi_b$ for the black open paths, the
link pattern $\pi_w$ for the white open paths, and the overall number
of black and white cycles, $\ell$ (a fact previously conjectured by
Cohn and Propp, also unpublished).

Call $\Psi_{n;\pm}(\pi_b, \pi_w; \ell)$ the number of FPL $\phi$ with
given triplet $(\pi_b, \pi_w; \ell)$, extending the definition of
$\Psi_{n;\pm}(\pi)$ used in the body of this paper.
Wieland proves that
\be
\Psi_{n;\pm}(\pi_b, \pi_w; \ell)
=
\Psi_{n;\pm}(R \pi_b, R^{-1} \pi_w; \ell)
\ef;
\ee
which, neglecting the refinement over $\pi_w$ and $\ell$, reduces to
\be
\Psi_{n;\pm}(\pi)
=
\Psi_{n;\pm}(R \pi)
\ef.
\ee
He reaches this result through a bijection $G$, called
\emph{gyration}, between the configurations in the pertinent refined
classes. This bijection operates locally over the elementary
plaquettes of the square lattice, and has the special property of
deforming only locally open monochromatic paths over the graph,
keeping fixed the endpoints of the intersection between the whole path
and the plaquette.

More precisely, the procedure operates in two steps
(cfr.\ \cite[sec.~5]{wie}), through two bijections $H_+$ and $H_-$,
each step involving the plaquettes with a given parity.  While $G$ is
a bijection from $\fpl(n,\pm)$ to itself, the two maps $H_+$ and $H_-$
are bijections mapping $\fpl(n,+)$ to $\fpl(n,-)$ and vice versa.
Then $G$ is obtained through the composition $G = H_- H_+$, $H_+^2 =
H_-^2 = 1$, and thus $G^{-1}=H_+ H_-$. The strong version of the
Wieland statement reads
\begin{align}
\label{eq.564675}
\Psi_{n;+}(\pi_b, \pi_w; \ell)
&=
\Psi_{n;-}(\pi_b, R^{-1} \pi_w; \ell)
\ef;
\\
\label{eq.564675b}
\Psi_{n;-}(\pi_b, \pi_w; \ell)
&=
\Psi_{n;+}(R \pi_b, \pi_w; \ell)
\ef;
\end{align}
Again, neglecting the refinement over $\pi_w$ and $\ell$, this reduces to
\begin{align}
\label{eq.564675X}
\Psi_{n;+}(\pi)
&=
\Psi_{n;-}(\pi)
\ef;
\\
\label{eq.564675bX}
\Psi_{n;-}(\pi)
&=
\Psi_{n;+}(R \pi)
\ef.
\end{align}
As a corollary, combining (\ref{eq.564675X}) with the discrete
reflection symmetry of the square along a vertical axis, for $n$
even, or along a diagonal, for $n$ odd, we also get 
$\Psi_{n;+}(\pi) = \Psi_{n;+}(V\pi)$, with 
$(i,j) \in \pi \leftrightarrow (2n+1-i,2n+1-j) \in V\pi$, completing
the statement on the dihedral symmetry of the enumerations.

Here we review Wieland proof, in a broader setting more suitable to
the generalizations we aim to. We do this in two main steps: in a
first moment, we concentrate on a single map $H$, inverting the
boundary conditions; at a later stage, we analyse how the construction
of a pair of distinct bijections is fruitfully exploited.

We consider a connected graph $\gc=(V,E)$ (not necessarily planar, and
with no given embedding on a surface). We require all vertices to have
degree $4$ or $2$, call $V' \subseteq V$ the set of degree-2 vertices,
and $E' \subseteq E$ the set of edges adjacent to $V'$. The existence
of degree-2 vertices, and absence of degree-1 vertices, apparently
seems at difference with the setting described in Section
\ref{sec.statem}, where we have degree-1 vertices on the boundary and
degree-4 vertices inside the region $\Lambda$: we will see later how
this case is recovered.

We define the set of valid FPL configurations on this graph,
$\fpl(\gc)$, as the set of maps $\phi \in \{b,w\}^{E}$ (black and
white), satisfying the ice-type constraint at all degree-4 vertices,
i.e.\ such that $\deg_b(v) = \deg_w(v)$ for each $v \in V \ssm V'$.
We also define the subsets $\fpl(\gc; \tau)$, of valid FPL
configurations $\phi \in \fpl(\gc)$, whose restriction to $E'$ is the
string $\tau \in \{b,w\}^{E'}$.  Call $V''(\tau) \subseteq V'$ the set
of vertices $v$ with $\deg_b(v)=\deg_w(v)=1$, and $2n$ its cardinality
(it is easily seen that $\fpl(\gc; \tau) = \emptyset$ if $|V''|$ is
odd). Label the vertices of $V''$ with indices from 1 to $2n$.

In such a domain, a configuration $\phi$ is composed of
monochromatic cycles, visiting only vertices in $V \ssm V''$, and
monochromatic open paths, having vertices in $V \ssm V''$ as interior
points, and vertices in $V''$ as endpoints. Given the labeling over
$V''$, $\phi$ determines a ``black'' matching $\pi_b$ of these $2n$
points, through the black open paths, and a ``white'' matching
$\pi_w$, through the white open paths. It also determines numbers
$\ell_b$, $\ell_w$ of black and white cycles, and
$\ell=\ell_b+\ell_w$. So, in particular it determines a triplet 
$(\pi_b, \pi_w; \ell)$, and we have refined enumerations 
$\Psi_{\gc; \tau}(\pi_b, \pi_w; \ell)$ for the cardinalities of these
classes inside $\fpl(\gc; \tau)$.

For a generic graph $\gc$, the matchings $\pi$ run over the whole set
of matchings over $2n$ points, of cardinality $(2n-1)!!$, while, if
$\gc$ has a planar embedding with all the vertices of $V''$ adjacent
to the same face (say, the external one), then we can restrict to link
patterns $\pi \in \LP(n)$.

Now consider a partition of $E$ into a collection of disjoint
unoriented cycles, $\Gamma=\{ \gamma_i \}$.
We want to construct a map $H_{\Gamma}$, that sends each $\phi \in
\fpl(\gc;\tau)$ to a $\phi' \in \fpl(\gc;\bar{\tau})$, such that, for
each cycle $\gamma \in\Gamma$, the three following conditions are
satisfied:
\begin{description}
\item[degree condition:] for $v \in \gamma$, we have two edges within
  $\gamma$ adjacent to $v$. If $0 \leq k \leq 2$ of them are black in
  $\phi$, then $k$ of them are white in $\phi'$.
\item[connectivity condition:] for $v, v' \in \gamma$, $v$ and $v'$
  are connected on $\gamma$ by an open black path in $\phi$ iff 
  they are connected on $\gamma$ by an open black path in $\phi'$, and
  similarly for white;
\item[alternation condition:] for $v \in V''(\tau)$, a single $\gamma$
  contains both adjacent edges. If in $\phi$ these edges are exactly
  one black and one white, then in $\phi'$ the black edge becomes
  white, and the white edge becomes black.
\end{description}

\noindent
We discuss later under which conditions on $\gc$ and $\Gamma$ one or
more maps $H_{\Gamma}: \fpl(\gc, \tau) \to \fpl(\gc, \bar{\tau})$
exist, satisfying these constraints, and when they are bijections.  We
note immediately that, because of the symmetry of the constraints,
(among which, the fact that they depend on $\tau$ only through
$V''(\tau)$, and $V''(\tau) = V''(\bar{\tau})$), if $H_{\Gamma}(\phi)$
is valid, over the domain $\fpl(\gc, \tau)$, the map
$\overline{H_{\Gamma}(\bar{\phi})}$ is valid on the domain $\fpl(\gc,
\bar{\tau})$, and, if $H_{\Gamma}(\phi)$ is a bijection, also
$\overline{H_{\Gamma}(\bar{\phi})}$ is a bijection.

Assuming that we have such a map $H_{\Gamma}$, we have
\begin{proposition}
\label{prop.HGam}
The FPL $\phi \in \fpl(\gc, \tau)$ and $\phi'=H_{\Gamma}(\phi) \in
\fpl(\gc, \bar{\tau})$ have the same triplet $(\pi_b, \pi_w; \ell)$.
\end{proposition}

\noindent
{\it Proof:} Remark that the degree condition ensures that $\phi' \in
\fpl(\gc)$, and the alternation condition ensures also that $\phi' \in
\fpl(\gc, \bar{\tau})$.  Take a monochromatic open path (m.o.p.) $p$
determined by $\phi$, say a black one. It is composed of an open
concatenation of m.o.p.'s contained into cycles $\gamma_i$ of $\Gamma$
(in a sequence $(i_1,i_2, \ldots)$ that may allow for repetitions, if
not consecutive). By the connectivity condition, m.o.p.'s within
cycles are sent to m.o.p.'s within cycles, with the same endpoints,
which thus still concatenate, leading overall to a m.o.p.\ $p'$ with
the same endpoints. This proves that $\pi_b(\phi) = \pi_b(\phi')$, and
analogously $\pi_w(\phi) = \pi_w(\phi')$. Now consider monochromatic
cycles of $\phi$. Take a cycle $c$, say black. Either it concides with
a cycle of $\Gamma$, in which case, by the degree condition, on
$\phi'$ we will have that $c$ is a white cycle, or it is composed of a
closed concatenation of m.o.p.'s contained into cycles $\gamma_i$ of
$\Gamma$ (at least two of them), in which case a reasoning analogous
to the one above for open paths allows to conclude that the endpoints
of the m.o.p.'s will be crossed in $\phi'$ by a monochromatic black
cycle $c'$. This proves that $\ell$ cannot decrease. But the three
conditions are symmetric w.r.t.\ $\phi$ and $\phi'$, so the reasonings
above can be repeated verbatim, starting with a cycle $c'$ in $\phi'$,
leading to the conclusion that $\ell$ cannot increase, thus it is
conserved.  \hfill $\square$

\vspace{2mm} 
\noindent
We now analyse how and when one can construct any map $H_{\Gamma}$
satisfying the three conditions above.  Remark that the conditions are
factorized over the cycles $\gamma_i$, so it suffices to concentrate
on single cycles of length $\ell$, and this makes feasible an analysis
for all graphs~$\gc$.

It turns out that the conditions can be satisfied if and only if all
the cycles $\gamma_i \in \Gamma$ have length $\ell$ at most $4$ (i.e.,
in the range $\{1,2,3,4\}$, as also loops are allowed), and all the
cycles adjacent to vertices in $V''$ have length at most $3$.

In these cases, the solution is unique, and leads to bijections,
except for the cycles of length 2 that are not adjacent to vertices in
$V'$, for which there are four solutions, two bijective and two
non-bijective ones.  We do not explore the most general case, but
limit ourselves to the most interesting one, in which, in this
ambiguous case, we take the single solution that satisfies also the
alternation condition (which is bijective). This ambiguity does not
appear in our domains $\Lambda$, as in this case we will not have
cycles of length $2$ not adjacent to the border.

The resulting solutions on the single cycle are involutive maps, that
just swap the black/white occupations of the edges, with the unique
exception of $\ell=4$, and the four edges in $\gamma$ having
alternating colouring, $(b,w,b,w)$ in cyclic order, in which case the
map acts as the identity (this is the only possibility in order to
preserve the connectivity constraint).

The fact that for $\ell>4$ there are no solutions is easily
proven. For example, take a configuration $\phi$ whose edges in
$\gamma$ are respectively
\[
(b,w,b,\underbrace{w,w,\ldots,w}_{\ell-3})
\ef.
\]
Then, the degree condition on the $\ell-4$ vertices internal to the
white path of length $\ell-3$ forces, for $\phi'$,
\[
(?,?,?,\underbrace{b,b,\ldots,b}_{\ell-3})
\ef,
\]
but this already breaks the connectivity condition, as we have two
points that are connected by a black arc in $\phi'$, but are not
in~$\phi$.

So, in conclusion, we have a precise set of conditions for the
existence of a bijection $H_{\Gamma}$, and a precise construction of
this bijection, that we summarize in a definition.
\begin{definition}
\label{def.triplet}
A triplet $(\gc, \tau, \Gamma)$ is \emph{valid}
if the following conditions are satisfied
\begin{itemize}
\item All the cycles of $\Gamma$ have length at most 4;
\item All the cycles of $\Gamma$ adjacent to a vertex in $V''(\tau)$
  have length at most 3.
\end{itemize}
In this positive case, the map $H_{\Gamma}: \fpl(\Lambda, \tau) \to
\fpl(\Lambda, \bar{\tau})$ is defined as follows. Calling
$\phi'=H_{\Gamma}(\phi)$,
\be
\label{eq.defHphi}
\phi'|_{\gamma_i}
=
\left\{
\begin{array}{ll}
\phi|_{\gamma_i}  
&  |\gamma_i|=4; \quad
\phi|_{\gamma_i} \in \big\{ (b,w,b,w), (w,b,w,b) \big\}; \\
\rule{0pt}{12pt}%
\overline{\phi}|_{\gamma_i}  & \textrm{otherwise.}
\end{array}
\right.
\ee
This very same definition also induces a map from
$\fpl(\Lambda, \bar{\tau})$ to $\fpl(\Lambda, \tau)$.
\end{definition}
Now we explain how to recover the original Wieland maps $H_{\pm}$ from
the map $H_{\Gamma}$ described above. Take the original FPL square
domain of side $n$. We thus have $4n$ vertices of degree 1 and $n^2$
vertices of degree 4. Label the degree-1 vertices from $1$ to $4n$ in
cyclic order. Colour the edges adjacent to these vertices according to
the alternating boundary condition~$\tau_+$.

Join together pairs of consecutive degree-1 vertices, i.e., for $H_+$,
glue together $(1,2)$, $(3,4)$, and so on, and, for $H_-$, glue
together $(4n,1)$, $(2,3)$, and so on. Now we have a graph $\gc$ with the
required properties, $|V'|=|V''|=2n$ and $|V \ssm V'|=n^2$.

It is easy to see that, in this graph, there is a single
possible choice of $\Gamma$, the one corresponding to take the
elementary plaquettes w.r.t.\ the planar embedding, of the given
parity that matches, on the boundary, the constraint that points in
$V'$ are covered (the unicity is indeed proven by starting from the
choices at the boundary, that are constrained, and continuing
recursively).

It is also easy to see that if we did not pair consecutive degree-1
vertices, w.r.t.\ the natural cyclic ordering and up to the trivial
ambiguity at the corners, there would have been no valid choice of
$\Gamma$. So, overall for all possible pairings of degree-1 vertices
in our square domain, there are only two valid glueing procedures and
associated $H_{\Gamma}$ maps, coinciding with~$H_{\pm}$.

The ambiguity at the corners is not totally negligible, and will have
a role in the following. It is however easier to visualize it as an
ambiguity in the planar embedding of the drawing, and then assume
that, for a given embedding, we only construct the two maps $H_{\pm}$
corresponding to the pairings along the cyclic ordering.

We now take our original square domain $\Lambda$, with boundary
condition $\tau_+$, glue together pairs of degree-1 vertices as
described above, apply $H_{\pm}$, and then split the degree-2 vertices
to recover the original domain $\Lambda$. The alternation condition
has now forced boundary conditions $\tau_-$. Furthermore, the
splitting of the vertices has caused an important ``switch'' among
black and white endpoints, w.r.t.\ their position in the $\Lambda$
domain: we started gluing together the $i$-th endpoint, black, and the
$(i\pm 1)$-th, white, we applied $H_{\pm}$, that satisfies the
alternation condition on the degree-2 vertex, so we end up splitting
the degree-2 vertex into the the $i$-th endpoint, now white, and the
$(i\pm 1)$-th one, now black.  This is what is responsible for the
rotations of the link patterns under 
$H_{\pm} : \fpl(n; +) \to \fpl(n; -)$, arising in (\ref{eq.564675})
and (\ref{eq.564675b}).

\subsection{Arbitrary regions $\Lambda$ and boundary
  conditions $\tau$: existence of the map}
\label{sec.LTexistence}

\noindent
Here we investigate regions $\Lambda$ which are portions of the
square lattice, take the pairing $(1,2)$, $(3,4)$, \ldots, of the
endpoints, producing a graph $\gc_+(\Lambda)$, and the set $\Gamma$
corresponding to the plaquettes of the lattice, with the appropriate
parity for covering the boundary points in $V'(\gc)$.

We already know from the definition \ref{def.triplet} under which
conditions a generic triplet $(\gc, \tau, \Gamma)$ is valid, as a set
of constraints on the length of the cycles in $\Gamma$. We want to
translate this to more effective conditions, in the special case of
the triplet $(\gc_+(\Lambda), \tau, \Gamma)$.

As all the plaquettes in the square lattice have length 4, the only
possibility for the triplet to be not valid is that we form long
cycles, adjacent to the border, in the glueing procedure. So we
concentrate on the boundary of $\Lambda$.

Consider the oriented boundary of $\Lambda$, $\partial \Lambda$, as a
closed path, say surrounding $Lambda$ in counter-clockwise
orientation, encoded as a sequence $\bm{\sigma}$ of ``steps'' in the
alphabet $\{-1,0,+1\}$, where 0 correspond to go straight, $+1$ to
rotate left (forming a \emph{convex} vertex in the polygon $\Lambda$),
and $-1$ to rotate right (forming a \emph{concave} vertex). A
necessary condition for $\partial \Lambda$ to be a closed path is that
this string has four $+1$ more than $-1$, and, for a rectangle,
$\bm{\sigma}$ is just composed of four $+1$ and some zeroes.

This string determines a sequence of $|\partial \Lambda|$ terminations
on the boundary of $\Lambda$, i.e.\ the edges in the set
$E_1(\Lambda)$ introduced in Section \ref{sec.defFPL}.  Terminations
interlace with the step, so that, say, the $k$-th step is between the
$k$-th termination and the $(k+1)$-th termination.  We have a simple
bijection of configurations, preserving the link pattern, if we
interchange the $k$-th termination and the $(k+1)$-th termination,
when the $k$-th step is a ``$+1$'', i.e. a convex corner, and this
possibility will be exploited in the following.

We must glue together terminations $2j-1$ and $2j$, so only the odd
values $\sigma_{2j-1}$ are relevant for the constraint.  In the
glueing of the pair above, if $\sigma_{2j-1}=0$ or $+1$, we form
cycles of length at most 3.  These cycles are always allowed.  If
instead $\sigma_{2j-1}=-1$ the glueing will form a cycle of length
4. This situation is interesting: if the two terminations have
different colour, then the vertex resulting from the glueing is in the
set $V''(\tau)$, and thus the triplet $(\gc(\Lambda),\tau,\Gamma)$ is
not valid. Conversely, if the two terminations have the same colour,
then the vertex is in the set $V' \ssm V''(\tau)$, and the triplet may
still be valid.

If, in $\bm{\sigma}$, we have consecutive $-1$ and $+1$ (in any
order), with the $-1$ in a position with odd index, the arbitrarity in
the ordering of the two terminations associated to $+1$ may be
critically exploited in order to produce a monochromatic pair glued
above the $-1$, and make the triplet valid.

\begin{figure}
\begin{align*}
&
\setlength{\unitlength}{15pt}
\begin{picture}(10,8)(0,0.5)
\put(0,0){\includegraphics[scale=3]{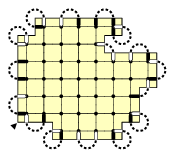}}
\put(6.9,6.1){\makebox[0pt][c]{$\scriptstyle{A}$}}
\put(9.5,4.1){\makebox[0pt][c]{$\scriptstyle{C}$}}
\put(8.6,3.2){\makebox[0pt][c]{$\scriptstyle{B}$}}
\end{picture}
&&
\setlength{\unitlength}{15pt}
\begin{picture}(10,8)(0,0.5)
\put(0,0){\includegraphics[scale=3]{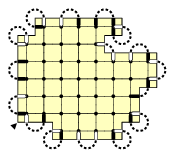}}
\put(6.9,6.1){\makebox[0pt][c]{$\scriptstyle{A}$}}
\put(9.5,4.1){\makebox[0pt][c]{$\scriptstyle{C}$}}
\put(8.6,3.2){\makebox[0pt][c]{$\scriptstyle{B}$}}
\end{picture}
\end{align*}

\setlength{\unitlength}{16pt}
\begin{picture}(21,3.5)(-2.2,-0.5)
\put(0,2){\makebox[0pt][c]{$\bm{\sigma}$}}
\put(1,2){\makebox[0pt][c]{$0$}}
\put(2,2){\makebox[0pt][c]{$-$}}
\put(3,2){\makebox[0pt][c]{$+$}}
\put(4,2){\makebox[0pt][c]{$0$}}
\put(5,2){\makebox[0pt][c]{$0$}}
\put(6,2){\makebox[0pt][c]{$0$}}
\put(7,2){\makebox[0pt][c]{$+$}}
\put(8,2){\makebox[0pt][c]{$-$}}
\put(9,2){\makebox[0pt][c]{$+$}}
\put(10,2){\makebox[0pt][c]{$0$}}
\put(11,2){\makebox[0pt][c]{$\stackrel{B}{-}$}}
\put(12,2){\makebox[0pt][c]{$\stackrel{C}{+}$}}
\put(13,2){\makebox[0pt][c]{$0$}}
\put(14,2){\makebox[0pt][c]{$+$}}
\put(15,2){\makebox[0pt][c]{$0$}}
\put(16,2){\makebox[0pt][c]{$0$}}
\put(17,2){\makebox[0pt][c]{$\stackrel{A}{-}$}}
\put(18,2){\makebox[0pt][c]{$-$}}
\put(19,2){\makebox[0pt][c]{$+$}}
\put(20,2){\makebox[0pt][c]{$+$}}
\put(21,2){\makebox[0pt][c]{$\cdots$}}

\put(0,1){\makebox[0pt][r]{left:}}
\put(1,1){\makebox[0pt][c]{$\underbrace{\circ \quad \bullet}$}}
\put(3,1){\makebox[0pt][c]{$\underbrace{\circ \quad \bullet}$}}
\put(5,1){\makebox[0pt][c]{$\underbrace{\circ \quad \circ}$}}
\put(7,1){\makebox[0pt][c]{$\underbrace{\bullet \quad \circ}$}}
\put(9,1){\makebox[0pt][c]{$\underbrace{\bullet \quad \circ}$}}
\put(11,1){\makebox[0pt][c]{$\underbrace{\bullet \quad \circ}$}}
\put(13,1){\makebox[0pt][c]{$\underbrace{\bullet \quad \circ}$}}
\put(15,1){\makebox[0pt][c]{$\underbrace{\bullet \quad \circ}$}}
\put(17,1){\makebox[0pt][c]{$\underbrace{\circ \quad \circ}$}}
\put(19,1){\makebox[0pt][c]{$\underbrace{\bullet \quad \circ}$}}

\put(0,0){\makebox[0pt][r]{right:}}
\put(1,0){\makebox[0pt][c]{$\underbrace{\circ \quad \bullet}$}}
\put(3,0){\makebox[0pt][c]{$\underbrace{\circ \quad \bullet}$}}
\put(5,0){\makebox[0pt][c]{$\underbrace{\circ \quad \circ}$}}
\put(7,0){\makebox[0pt][c]{$\underbrace{\bullet \quad \circ}$}}
\put(9,0){\makebox[0pt][c]{$\underbrace{\bullet \quad \circ}$}}
\put(11,0){\makebox[0pt][c]{$\underbrace{\bullet \quad \bullet}$}}
\put(13,0){\makebox[0pt][c]{$\underbrace{\circ \quad \circ}$}}
\put(15,0){\makebox[0pt][c]{$\underbrace{\bullet \quad \circ}$}}
\put(17,0){\makebox[0pt][c]{$\underbrace{\circ \quad \circ}$}}
\put(19,0){\makebox[0pt][c]{$\underbrace{\bullet \quad \circ}$}}
\end{picture}

\caption{\label{fig.exTriplet}On the left, an example of a triplet
  $(\gc_+(\Lambda), \tau, \Gamma)$ which is not valid. Indeed, while
  the arc $A$, over a concave corner, is not violating the condition,
  because the two terminations are monochromatic, the arc $B$ is
  violating the condition. However, this concave angle is adjacent to
  a non-monochromatic convex angle (next to letter $C$), so we can
  exploit the invariance of the system under swap of the two
  terminations over $C$. The resulting valid triplet $(\gc_+(\Lambda),
  \tau, \Gamma)$ is shown on the right. The string $\bm{\sigma}$, and
  the terminations for the domains on the left and right side of the
  picture, are also descibed in the table.
}
\end{figure}

These conditions are illustrated in an example in
figure~\ref{fig.exTriplet}.

\subsection{An extended diagram algebra}
\label{sec.TLca}

\noindent
In Section \ref{sec.TL} we introduce the representation of the $2n$
(affine) Temperley-Lieb operators $e_j$, $1 \leq j \leq 2n$, acting
over link patterns $\pi \in \LP(n)$. Here we introduce ``diagram''
operators $c_j$ and $a_j$ that relate spaces $\LP(n)$ with different
values of $n$. In order to make the analysis pictorially simple, we do
\emph{not} introduce the ``affine'' version of this algebra, and only
introduce operators $c_j$ with $1 \leq j \leq 2n-1$ and $a_j$ with $1
\leq j \leq 2n+1$.

We thus recall the diagram definitions of $R$ and of $e_j$
\begin{align}
R: \quad &\LP(n) \to \LP(n)
&&
\setlength{\unitlength}{10pt}
\raisebox{-10pt}{
\begin{picture}(11,2.5)
\put(0,1){\includegraphics[scale=2]{fig_TL_R.eps}}
\put(1,0){\makebox[0pt][c]{$\scriptstyle{1}$}}
\put(2,0){\makebox[0pt][c]{$\scriptstyle{2}$}}
\put(3,0){\makebox[0pt][c]{$\scriptstyle{3}$}}
\put(4,0){\makebox[0pt][c]{$\scriptstyle{\cdots}$}}
\put(10,0){\makebox[0pt][c]{$\scriptstyle{2n}$}}
\put(7,1.8){\makebox[0pt][c]{$\scriptstyle{\cdots}$}}
\end{picture}
}
\\
e_j: \quad &\LP(n) \to \LP(n)
&&
\setlength{\unitlength}{10pt}
\raisebox{-10pt}{
\begin{picture}(11,3.5)
\put(0,1){\includegraphics[scale=2]{fig_TL_ej.eps}}
\put(1,0){\makebox[0pt][c]{$\scriptstyle{1}$}}
\put(2,0){\makebox[0pt][c]{$\scriptstyle{2}$}}
\put(3,0){\makebox[0pt][c]{$\scriptstyle{3}$}}
\put(4,0){\makebox[0pt][c]{$\scriptstyle{\cdots}$}}
\put(5.7,0){\makebox[0pt][c]{$\scriptstyle{j}$}}
\put(7.2,0){\makebox[0pt][c]{$\scriptstyle{j+1}$}}
\put(9,0){\makebox[0pt][c]{$\scriptstyle{\cdots}$}}
\put(10,0){\makebox[0pt][c]{$\scriptstyle{2n}$}}
\put(4,1.5){\makebox[0pt][c]{$\scriptstyle{\cdots}$}}
\put(9,1.5){\makebox[0pt][c]{$\scriptstyle{\cdots}$}}
\end{picture}
}
\intertext{We now give the definitions of $c_j$ and $a_j$ (\emph{close} and
\emph{add})}
c_j: \quad &\LP(n) \to \LP(n-1)
&&
\setlength{\unitlength}{10pt}
\raisebox{-10pt}{
\begin{picture}(11,3)
\put(0,1){\includegraphics[scale=2]{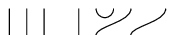}}
\put(1,0){\makebox[0pt][c]{$\scriptstyle{1}$}}
\put(2,0){\makebox[0pt][c]{$\scriptstyle{2}$}}
\put(3,0){\makebox[0pt][c]{$\scriptstyle{3}$}}
\put(4,0){\makebox[0pt][c]{$\scriptstyle{\cdots}$}}
\put(5.8,2.9){\makebox[0pt][c]{$\scriptstyle{j}$}}
\put(7.,2.9){\makebox[0pt][c]{$\scriptstyle{j+1}$}}
\put(8.4,2.9){\makebox[0pt][c]{$\scriptstyle{j+2}$}}
\put(7,0){\makebox[0pt][c]{$\scriptstyle{\cdots}$}}
\put(8.6,0){\makebox[0pt][c]{$\scriptstyle{2n-2}$}}
\put(10,2.9){\makebox[0pt][c]{$\scriptstyle{2n}$}}
\put(4,1.5){\makebox[0pt][c]{$\scriptstyle{\cdots}$}}
\put(8,1.5){\makebox[0pt][c]{$\scriptstyle{\cdots}$}}
\put(6,0){\makebox[0pt][c]{$\scriptstyle{j}$}}
\end{picture}
}
\\
a_j: \quad &\LP(n) \to \LP(n+1)
&&
\setlength{\unitlength}{10pt}
\raisebox{-10pt}{
\begin{picture}(11,4)
\put(0,1){\includegraphics[scale=2]{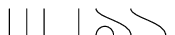}}
\put(1,0){\makebox[0pt][c]{$\scriptstyle{1}$}}
\put(2,0){\makebox[0pt][c]{$\scriptstyle{2}$}}
\put(3,0){\makebox[0pt][c]{$\scriptstyle{3}$}}
\put(4,0){\makebox[0pt][c]{$\scriptstyle{\cdots}$}}
\put(5.8,0){\makebox[0pt][c]{$\scriptstyle{j}$}}
\put(7.,0){\makebox[0pt][c]{$\scriptstyle{j+1}$}}
\put(8.4,0){\makebox[0pt][c]{$\scriptstyle{j+2}$}}
\put(7,2.9){\makebox[0pt][c]{$\scriptstyle{\cdots}$}}
\put(8.6,2.9){\makebox[0pt][c]{$\scriptstyle{2n-2}$}}
\put(10,0){\makebox[0pt][c]{$\scriptstyle{2n}$}}
\put(4,1.5){\makebox[0pt][c]{$\scriptstyle{\cdots}$}}
\put(8.35,1.5){\makebox[0pt][c]{$\scriptstyle{\cdots}$}}
\put(6,2.9){\makebox[0pt][c]{$\scriptstyle{j}$}}
\end{picture}
}
\end{align}
We have a number of algebraic relations, easily deduced from the
drawing of these diagrams, in a fashion similar to the deduction of
the Temperley-Lieb relations (\ref{eqs.TL}). We do not list all of
them, but only a subset that will be used in the following. We clearly
have
\begin{align}
\label{eq.cae}
c_j a_j 
&= 1
\ef;
\\
\label{eq.ace}
a_j c_j 
&= e_j
\ef;
\\
[e_j,c_k]
&=
[e_j,a_k] = 0
&&
\textrm{if $k-j \geq 2$.}
\end{align}
From which we get, for a set $J \subset \{ 1,2,\ldots, 2n-1 \}$ with
no pairs of consecutive indices,
\be
\label{eq.AACCE}
\Big(
\prod_{\substack{j \in J \\ \textrm{decreasing}\\ \textrm{order}}}
a_j
\Big)
\Big(
\prod_{\substack{j \in J \\ \textrm{increasing}\\ \textrm{order}}}
c_j
\Big)
=
\prod_{j \in J}
e_j
\ef.
\ee
(remark that the product on the right-hand side is well defined
without ordering prescriptions, as the absence of consecutive pairs
implies that the operators $e_j$ in the monomial do commute).

An example of application of these operators has already been
anticipated in the simple path reversal relation,
Proposition~\ref{prop.SPR}, where we had
\be
\label{eq.postprop.SPR}
e_{2 m-1} \ket{s_1} = a_{2m-1} \ket{s_2}
\ef.
\ee
Similarly, using 
(\ref{eq.cae}) and (\ref{eq.ace}),
we could also write a relation in $\mathbb{C}^{\LP(m-1)}$
\be
\label{eq.postprop.SPR.2}
c_{2 m-1} \ket{s_1} = \ket{s_2}
\ef.
\ee

\subsection{Gyration in arbitrary regions $\Lambda$ and boundary
  conditions $\tau$}

\noindent
Here we analyse the consequences of Proposition \ref{prop.HGam} on the
refined enumerations of FPL on the ensembles $\fpl(\Lambda,\tau)$ and
$\fpl(\Lambda,\bar{\tau})$, with $\Lambda$ as in Section
\ref{sec.statem}, and $\Gamma$ induced by the pairing $(1,2)$,
$(3,4)$, \ldots, of consecutive endpoints in cyclic order.  As we will
see, the notation introduced in Section \ref{ssec.VectNot} is
specially suitable at this purpose.

We will assume here that the triplet $(\gc(\Lambda),\tau,\Gamma)$ is
valid, a question already addressed in Section~\ref{sec.LTexistence}.

The points on the boundary are collected into $N$ pairs
$\{(2i-1,2i)\}_{1 \leq i \leq N}$, of which $N_{bb}$, $N_{ww}$,
$N_{bw}$ and $N_{wb}$ ones coloured, in cyclic order, as $(b,b)$,
$(w,w)$, $(b,w)$ and $(w,b)$ respectively.  The analysis of the
gyration operation will produce a relation among vectors in the linear
space $\mathbb{C}^{\LP(n)}$, with $n=N_{bw}+N_{wb}$.  In the domains
$(\Lambda,\tau)$ and $(\Lambda,\bar{\tau})$ we have respectively $2
m_1 = 2N_{bb} + N_{bw}+N_{wb}$ and $2 m_2 = 2N_{ww} + N_{bw}+N_{wb}$
black terminations, with $m_1, m_2 \geq n$. So, we have enumerations
$\Psi_{\Lambda,\tau}(\pi_1)$ and $\Psi_{\Lambda,\bar{\tau}}(\pi_2)$
with $\pi_1 \in \LP(m_1)$ and $\pi_2 \in \LP(m_2)$.
We will define the two states in $\mathbb{C}^{\LP(m_1)}$ and
$\mathbb{C}^{\LP(m_2)}$
\begin{align}
\label{eq.defs1s2}
\ket{s_1}
&=
\sum_{\pi_1 \in \LP(m_1)}
\Psi_{\Lambda,\tau}(\pi_1)
\ket{\pi_1}
\ef;
&
\ket{s_2}
&=
\sum_{\pi_2 \in \LP(m_2)}
\Psi_{\Lambda,\bar{\tau}}(\pi_2)
\ket{\pi_2}
\ef.
\end{align}
Call $J_1$ and $J_2$ the following sets of indices, with cardinalities
$|J_1|=N_{bb}$ and $|J_2|=N_{ww}$, in the sets
$\{1,\ldots,\hbox{$2m_1-1$} \}$ and $\{1,\ldots, \hbox{$2m_2-1$} \}$.
The set $J_1$ collects, for each pair $(2i-1,2i)$ that is coloured
$(b,b)$, the index of the left-most termination, according to the
cyclic labeling of the $2m_1$ black terminations, in the domain
$(\Lambda, \tau)$. The set $J_2$ does the analogous thing, for the
domain $(\Lambda, \bar{\tau})$. Remark that $J_{1,2}$ do not contain
pairs of consecutive indices, as a left-most termination of a
monochromatic pair, is followed by the right-most termination of the
same pair, that thus is not in $J_{1,2}$. Also remark that $2m_1
\not\in J_1$, as the last termination is either in a $(b,w)$ or a
$(w,b)$ pair, or it is the right-most termination of a $(b,b)$ pair.

The consequence of the general procedure described in the previous
section reads
\begin{proposition}[Generalized gyration]
\label{prop.gengyr}
\be
\label{eq.prop.gengyr}
\Big(
\prod_{\substack{j \in J_1 \\ \mathrm{increasing}\\ \mathrm{order}}}
c_j
\Big)
\ket{s_1}
=
\Big(
\prod_{\substack{j \in J_2 \\ \mathrm{increasing}\\ \mathrm{order}}}
c_j
\Big)
\ket{s_2}
\ef.
\ee
\end{proposition}
This equation should be considered as the generalization of equation
(\ref{eq.564675X}), that in vector notation, and with the definitions
(\ref{eq.defs1s2}), just reads
$\ket{s_1} = \ket{s_2}$. Indeed, in the simpler case of
$N_{bb}=N_{ww}=0$, we have $m_1=m_2=n$,
and $J_1 = J_2 = \emptyset$.

\vspace{2mm}
\noindent
{\it Proof of Proposition \ref{prop.gengyr}.}  We should relate the
domains $(\Lambda, \tau)$ and $(\Lambda, \bar{\tau})$ to the
modifications $(\gc, \tau)$ and $(\gc, \bar{\tau})$ in which the
endpoints are glued pairwise.  Then, from the assumption that the
triplet $(\gc(\Lambda),\tau,\Gamma)$ is valid, (which also implied
that $(\gc(\Lambda),\bar{\tau},\Gamma)$ is valid, by the observation
in Definition \ref{def.triplet}), we can apply the statement of
Proposition~\ref{prop.HGam}.  This statement provides information only
for enumerations refined accordingly to link patterns $\pi \in
\LP(n)$, as we have $|V''|=2n$ on $\gc$.  Read back on $\Lambda$, it
provides information only on connectivity properties of the black
endpoints in pairs $(b,w)$ or $(w,b)$, provided that the endpoints in
pairs $(b,b)$ are glued together (similarly, it provides information
on connectivity properties of the white endpoints in pairs $(b,w)$ or
$(w,b)$, provided that the endpoints in pairs $(w,w)$ are glued
together, and also information on the number $\ell$ of cycles, but we
do not use this information here).  Glueing together a pair $(b,b)$
corresponds exactly to apply an operator $c_j$. The product is
performed putting $c_j$ factors with higher index more on the right
(they act before on the state), at the aim of reading easily the
labeling of the endpoints (indeed, for $k-j \geq 2$, we have $c_j c_k
= c_{k-2} c_j$, at difference with $e_j e_k = e_k e_j$).  \hfill
$\square$

The idea beyond the use of relation (\ref{eq.prop.gengyr}) for proving
the Razumov-Stroganov conjecture, that involves operators $e_j$, is to
exploit combinations of these relations, especially comparing the
result of the two gyration operations $H_{\pm}$, and then multiply
both sides of the resulting relation by an appropriate monomial in the
$a_j$'s, reproducing Temperley-Lieb operators as a consequence of
equation (\ref{eq.AACCE}).

This cannot be done in general (the appropriate monomial for each side
of the equation is unique, and may be different on the two sides), but
this will be the case for the special situations described by
Proposition \ref{prop.gyrat}, as we describe in
Section~\ref{sec.gyraproofs}.


\subsection{Proof of the gyration relations}
\label{sec.gyraproofs}

\noindent
Here we specialize the general statements proven in the previous
sections to the states defined in Section \ref{ssec.auxil}, in order
to prove the gyration relations collected in
Proposition~\ref{prop.gyrat}.

\begin{figure}[tb]
\[
\begin{array}{cc}
e_j \ssb{c}{j} 
&
e_{j-1} \ssb{c}{j} 
\\
\setlength{\unitlength}{10pt}
\begin{picture}(15,5)(-.5,0)
\put(0.25,0.5){\includegraphics[scale=2]{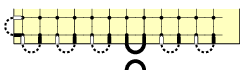}}
\put(0.2,2.3){\makebox[0pt][c]{$\scriptstyle{2n}$}}
\put(2.5,1){\makebox[0pt][c]{$\scriptstyle{1}$}}
\put(4.5,1){\makebox[0pt][c]{$\scriptstyle{2}$}}
\put(5.5,1){\makebox[0pt][c]{$\scriptstyle{\cdots}$}}
\put(6.5,1){\makebox[0pt][c]{$\scriptstyle{j-1}$}}
\put(7.3,0){\makebox[0pt][c]{$\scriptstyle{j}$}}
\put(8.7,0){\makebox[0pt][c]{$\scriptstyle{j+1}$}}
\put(10.5,1){\makebox[0pt][c]{$\scriptstyle{j+2}$}}
\put(11.9,1){\makebox[0pt][c]{$\scriptstyle{\cdots}$}}
\end{picture}
&
\setlength{\unitlength}{10pt}
\begin{picture}(15,5)(-.5,0)
\put(0.25,0.5){\includegraphics[scale=2]{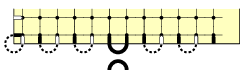}}
\put(0.2,2.8){\makebox[0pt][c]{$\scriptstyle{2n}$}}
\put(2.5,1){\makebox[0pt][c]{$\scriptstyle{1}$}}
\put(4.5,1){\makebox[0pt][c]{$\scriptstyle{2}$}}
\put(5.5,1){\makebox[0pt][c]{$\scriptstyle{\cdots}$}}
\put(6.3,0){\makebox[0pt][c]{$\scriptstyle{j-1}$}}
\put(7.7,0){\makebox[0pt][c]{$\scriptstyle{j}$}}
\put(8.5,1){\makebox[0pt][c]{$\scriptstyle{j+1}$}}
\put(10.5,1){\makebox[0pt][c]{$\scriptstyle{j+2}$}}
\put(11.9,1){\makebox[0pt][c]{$\scriptstyle{\cdots}$}}
\end{picture}
\\
\setlength{\unitlength}{10pt}
\begin{picture}(15,5)(-.5,0)
\put(0.25,0.5){\includegraphics[scale=2]{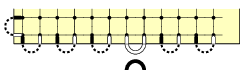}}
\put(0.2,3.8){\makebox[0pt][c]{$\scriptstyle{2n}$}}
\put(1.5,1){\makebox[0pt][c]{$\scriptstyle{1}$}}
\put(3.5,1){\makebox[0pt][c]{$\scriptstyle{2}$}}
\put(4.5,1){\makebox[0pt][c]{$\scriptstyle{\cdots}$}}
\put(5.5,1){\makebox[0pt][c]{$\scriptstyle{j-1}$}}
\put(7.3,0){\makebox[0pt][c]{$\scriptstyle{j}$}}
\put(8.7,0){\makebox[0pt][c]{$\scriptstyle{j+1}$}}
\put(9.5,1){\makebox[0pt][c]{$\scriptstyle{j+2}$}}
\put(10.9,1){\makebox[0pt][c]{$\scriptstyle{\cdots}$}}
\end{picture}
&
\setlength{\unitlength}{10pt}
\begin{picture}(15,5)(-.5,0)
\put(0.25,0.5){\includegraphics[scale=2]{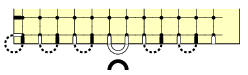}}
\put(0,3.5){\makebox[0pt][c]{$\scriptstyle{2n-1}$}}
\put(1.5,1){\makebox[0pt][c]{$\scriptstyle{2n}$}}
\put(3.5,1){\makebox[0pt][c]{$\scriptstyle{1}$}}
\put(4.5,1){\makebox[0pt][c]{$\scriptstyle{\cdots}$}}
\put(6.3,0){\makebox[0pt][c]{$\scriptstyle{j-1}$}}
\put(7.7,0){\makebox[0pt][c]{$\scriptstyle{j}$}}
\put(5.5,1){\makebox[0pt][c]{$\scriptstyle{j-2}$}}
\put(9.5,1){\makebox[0pt][c]{$\scriptstyle{j+1}$}}
\put(10.9,1){\makebox[0pt][c]{$\scriptstyle{\cdots}$}}
\end{picture}
\end{array}
\]
\caption{Top: the states $e_j \ssb{c}{j}$ and $e_{j-1} \ssb{c}{j}$. As
  the monochromatic pairs are glued together by the Temperley-Lieb
  operators, we can apply the generalized gyration, $H_+$ on 
  $e_j \ssb{c}{j}$ and $H_-$ on $e_{j-1} \ssb{c}{j}$. The result is
  shown on the bottom. The resulting domain are identical, up to an overall
  rotation of the indices.\label{fig.gyraC}}
\end{figure}

We start with equation (\ref{eq.2354765a}), concerning a state
$\ssb{c}{j}$. We analyse it in its ``frozen'' version, on a
rectangular $(n-1) \times n$ domain, with alternating boundary
conditions except for three consecutive black terminations, with
indices $j-1$, $j$ and $j+1$. As we have no concave angles, the
triplets $(\gc(\Lambda),\tau,\Gamma_{\pm})$ are automatically valid
for both pairings.  Also, for both pairings $\Gamma_{\pm}$ we have
$N_{bb}=1$ and $N_{ww}=0$, thus $|J_{1,\pm}|=1$ and $|J_{2,\pm}|=0$.
A direct inspection shows that $J_{1,+}=\{j-1\}$ and $J_{2,+}=\{j\}$,
and, combined together the result of Proposition \ref{prop.gengyr} in
the two cases, we get
\be
\label{eq.76565486}
c_j \ssb{c}{j}
=
R c_{j-1} \ssb{c}{j}
\ef.
\ee
Multiplying both sides by $a_j$ we get
\be
a_j c_j \ssb{c}{j}
=
a_j R c_{j-1} \ssb{c}{j}
=
R a_{j-1} c_{j-1} \ssb{c}{j}
\ef.
\ee
(Strictly speaking, the last passage can be done only for $j>2$. For
$j=1$, we could solve the apparent problem by performing a rotation of
the indices at the beginning.)
Thus, using (\ref{eq.ace}),
\be
e_j \ssb{c}{j}
=
R e_{j-1} \ssb{c}{j}
\ef,
\ee
as was to be proven.
These steps are illustrated in figure~\ref{fig.gyraC}.

We now analyse equation (\ref{eq.2354765aa}), concerning a state
$\ssb{a}{j}$. Again, we analyse the frozen domain, which now is the $n
\times n$ square, with part of the bottom row removed (namely all the
sites on the right of the decimated one), This domain has a single
concave turning in the perimeter, adjacent to a convex turning (i.e.,
a substring $\ldots,0,0,+1,-1,0,0,\ldots$ in $\bm{\sigma}$). The
boundary conditions are alternating, up to possibly using the
ambiguity in the ordering of the terminations at convex angles.

Indeed, as we have a concave angle, in one of the two maps $H_{\pm}$
(precisely, in $H_+$) we \emph{need} to exploit the ambiguity, before
performing the map, and
swap the terminations at the convex corner adjacent to the concave
one, in order to produce a valid triplet
$(\gc(\Lambda),\tau,\Gamma_H)$. So, in this case we have
$N_{bb}=N_{ww}=1$, and $J_1=J_2=\{j\}$.

\begin{figure}[tb]
\[
\begin{array}{cc}
\setlength{\unitlength}{10pt}
\begin{picture}(15,5)(-.5,0)
\put(-0.25,0.5){\includegraphics[scale=2]{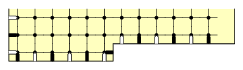}}
\put(-0.2,2.8){\makebox[0pt][c]{$\scriptstyle{2n}$}}
\put(1,0){\makebox[0pt][c]{$\scriptstyle{1}$}}
\put(3,0){\makebox[0pt][c]{$\scriptstyle{2}$}}
\put(3.8,0){\makebox[0pt][c]{$\scriptstyle{\cdots}$}}
\put(5,0){\makebox[0pt][c]{$\scriptstyle{j-1}$}}
\put(6.8,1){\makebox[0pt][c]{$\scriptstyle{j}$}}
\put(8.2,1){\makebox[0pt][c]{$\scriptstyle{j+1}$}}
\put(10,1){\makebox[0pt][c]{$\scriptstyle{j+2}$}}
\put(11.5,1){\makebox[0pt][c]{$\scriptstyle{\cdots}$}}
\end{picture}
&
\setlength{\unitlength}{10pt}
\begin{picture}(15,5)(-.5,0)
\put(-0.25,0.5){\includegraphics[scale=2]{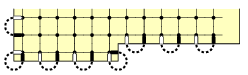}}
\put(-0.2,2.8){\makebox[0pt][c]{$\scriptstyle{2n}$}}
\put(1,0){\makebox[0pt][c]{$\scriptstyle{1}$}}
\put(3,0){\makebox[0pt][c]{$\scriptstyle{2}$}}
\put(3.8,0){\makebox[0pt][c]{$\scriptstyle{\cdots}$}}
\put(5,0){\makebox[0pt][c]{$\scriptstyle{j-1}$}}
\put(6.8,1){\makebox[0pt][c]{$\scriptstyle{j}$}}
\put(8.2,1){\makebox[0pt][c]{$\scriptstyle{j+1}$}}
\put(10,1){\makebox[0pt][c]{$\scriptstyle{j+2}$}}
\put(11.5,1){\makebox[0pt][c]{$\scriptstyle{\cdots}$}}
\end{picture}
\\
\setlength{\unitlength}{10pt}
\begin{picture}(15,5)(-.5,0)
\put(-0.25,0.){\includegraphics[scale=2]{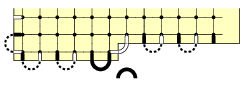}}
\put(-0.2,2.8){\makebox[0pt][c]{$\scriptstyle{2n}$}}
\put(1,0){\makebox[0pt][c]{$\scriptstyle{1}$}}
\put(3,0){\makebox[0pt][c]{$\scriptstyle{2}$}}
\put(3.8,0){\makebox[0pt][c]{$\scriptstyle{\cdots}$}}
\put(6.2,-0.5){\makebox[0pt][c]{$\scriptstyle{j-1}$}}
\put(7.7,-0.5){\makebox[0pt][c]{$\scriptstyle{j}$}}
\put(8.2,1){\makebox[0pt][c]{$\scriptstyle{j+1}$}}
\put(10,1){\makebox[0pt][c]{$\scriptstyle{j+2}$}}
\put(11.5,1){\makebox[0pt][c]{$\scriptstyle{\cdots}$}}
\end{picture}
&
\setlength{\unitlength}{10pt}
\begin{picture}(15,5)(-.5,0)
\put(-0.25,0.5){\includegraphics[scale=2]{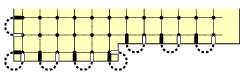}}
\put(-0.2,3.8){\makebox[0pt][c]{$\scriptstyle{2n}$}}
\put(-0.2,1.8){\makebox[0pt][c]{$\scriptstyle{1}$}}
\put(2,0){\makebox[0pt][c]{$\scriptstyle{2}$}}
\put(2.8,0){\makebox[0pt][c]{$\scriptstyle{\cdots}$}}
\put(4,0){\makebox[0pt][c]{$\scriptstyle{j-1}$}}
\put(6,0){\makebox[0pt][c]{$\scriptstyle{j}$}}
\put(7.2,1){\makebox[0pt][c]{$\scriptstyle{j+1}$}}
\put(9,1){\makebox[0pt][c]{$\scriptstyle{j+2}$}}
\put(10.5,1){\makebox[0pt][c]{$\scriptstyle{\cdots}$}}
\end{picture}
\\
\setlength{\unitlength}{10pt}
\begin{picture}(15,5)(-.5,0)
\put(-0.25,0.){\includegraphics[scale=2]{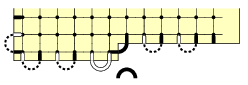}}
\put(-0.6,3.8){\makebox[0pt][c]{$\scriptstyle{2n-1}$}}
\put(-0.2,1.8){\makebox[0pt][c]{$\scriptstyle{2n}$}}
\put(2,0){\makebox[0pt][c]{$\scriptstyle{1}$}}
\put(4,0){\makebox[0pt][c]{$\scriptstyle{2}$}}
\put(4.8,0){\makebox[0pt][c]{$\scriptstyle{\cdots}$}}
\put(6.2,-0.5){\makebox[0pt][c]{$\scriptstyle{j-1}$}}
\put(7.7,-0.5){\makebox[0pt][c]{$\scriptstyle{j}$}}
\put(9.2,1){\makebox[0pt][c]{$\scriptstyle{j+1}$}}
\put(11,1){\makebox[0pt][c]{$\scriptstyle{j+2}$}}
\put(12.5,1){\makebox[0pt][c]{$\scriptstyle{\cdots}$}}
\end{picture}
&
\setlength{\unitlength}{10pt}
\begin{picture}(15,5)(-.5,0)
\put(-0.25,0.5){\includegraphics[scale=2]{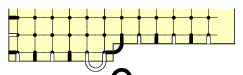}}
\put(-0.2,3.8){\makebox[0pt][c]{$\scriptstyle{2n}$}}
\put(-0.2,1.8){\makebox[0pt][c]{$\scriptstyle{1}$}}
\put(2,0){\makebox[0pt][c]{$\scriptstyle{2}$}}
\put(2.8,0){\makebox[0pt][c]{$\scriptstyle{\cdots}$}}
\put(4,0){\makebox[0pt][c]{$\scriptstyle{j-1}$}}
\put(6.2,-0.5){\makebox[0pt][c]{$\scriptstyle{j}$}}
\put(7.7,-0.5){\makebox[0pt][c]{$\scriptstyle{j+1}$}}
\put(9,1){\makebox[0pt][c]{$\scriptstyle{j+2}$}}
\put(10.5,1){\makebox[0pt][c]{$\scriptstyle{\cdots}$}}
\end{picture}
\end{array}
\]
\caption{Two different manipulations of the state $\ssb{a}{j}$. Left
  column: before applying $H_+$ we need to swap the two terminations
  on the convex angle adjacent to the concave one. Then, we can glue
  the terminations in monochromatic pairs (a black and a white one),
  and perform gyration. Right column: we can perform gyration
  immediately. In order to compare with the resulting domain on the
  left column, we have the possibility of swapping the two
  terminations, and glueing the adjacent monochromatic pairs. A
  rotation overall has resulted.\label{fig.gyraA}}
\end{figure}

We do not need to swap the terminations when applying the other map,
$H_-$. However, at the aim of combining the result of relation
(\ref{eq.prop.gengyr}) for the two maps $H_{\pm}$, and have identical
domains, we need to swap the terminations \emph{after} the application
of $H_-$, and then apply an operator $c_j$ on both sides of the
relation. This leads to the equation
\be
c_j \ssb{a}{j}
=
R c_{j-1} \ssb{a}{j}
\ef.
\ee
For reasonings identical to the ones following equation
(\ref{eq.76565486}), we thus get
\be
e_j \ssb{a}{j}
=
R e_{j-1} \ssb{a}{j}
\ef,
\ee
as was to be proven.
These steps are illustrated in figure~\ref{fig.gyraA}

The proof for the other three equations in Proposition
\ref{prop.gyrat} is very similar to the one for equation
(\ref{eq.2354765aa}), and we omit it.
\hfill $\square$

\section{Perspectives of generalization}
\label{sec.persp}

\noindent
In this paper, from the very beginning in Section \ref{sec.defFPL}, we
defined FPL configurations on portions of the square lattice. However,
the reader may have noticed that there is much space for
generalizations. This is clearly the case for our approach to
gyration, in Section \ref{sec.wie}, and also for a crucial step of the
proof, constituted by Lemma \ref{thm.HVgen} and Corollary~\ref{cor.HV}.

\begin{figure}
\begin{center}
\includegraphics[scale=2.7]{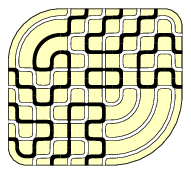}
\includegraphics[scale=2.5]{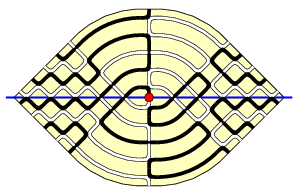}
\end{center}
\caption{Left: a FPL in a domain which is not a portion of the square
  lattice (because it has four triangular plaquettes), but shows both
  dihedral symmetry in the refined enumeration, and Razumov-Stroganov
  correspondence of these enumerations with the ground-state
  wavefunction of the periodic $O(1)$ loop model.  Right: a FPL in a domain
  which is not a portion of the square lattice (because it has an
  internal vertex with degree 2). This configuration is symmetric
  under reflection w.r.t.\ the horizontal axis, combined with
  complementation. The enumerations of the symmetric configurations
  show Razumov-Stroganov correspondence with the ground-state
  wavefunction of the open $O(1)$ loop model.\label{fig.multicone}}
\end{figure}

Consider for example an ensemble of FPL illustrated by the
configuration in figure~\ref{fig.multicone}, left.
It is the case that this domain has FPL enumerations with dihedral
symmetry, and also that these enumerations are proportional to the
very same set of integers for the square of side $n=10$, and thus,
also for the $O(1)$ loop model with $2n$ sites (however, note that
there is a non-trivial integer proportionality factor).

It would not be hard to show that the very same line of proof in this
paper works for proving that $(H_n - 2n) \ket{s_{\Lambda}}=0$ also for
this domain $\Lambda$. However, we postpone this analysis to a
different paper \cite{usprep}, where we also undertake the more
ambitious task of exausting the classification of the possible
structures for which the gyration mechanism works, i.e.\ all the
graphs for which the FPL refined enumerations have dihedral symmetry.

A similar goal will be accomplished also for further refinements of
FPL, as conjectured in \cite{Raz-Str-2} for the case of Half-Turn and
Quarter-Turn Symmetric FPL, and for Vertically-Symmetric FPL on the
square, for which it has been noticed that the link-pattern
enumerations are related to the integers in the ground-state
wavefunction for the closed or open system \cite{Raz-Str-2} (instead
of the periodic system, the case at hand in this paper).  Indeed, the
broader family of domains depicted above, all showing dihedral
symmetry, and the Razumov-Stroganov correspondence with the periodic
$O(1)$ loop model, contains subfamilies with an involutive symmetry
(it may be a reflection, or a rotation by 180 degrees, possibly
combined with a complementation). We will show how, similarly to how
the domain on the left part of figure~\ref{fig.multicone} generalizes
the $n \times n$ square domain, these subfamilies with involutive
symmetries generalize the Half-Turn symmetric, Quarter-Turn symmetric,
Vertically-Symmetric FPL domains, and the enumerations of symmetric
FPL within the symmetric domains see the emergence of
Razumov-Stroganov correspondence with the \emph{closed} or \emph{open}
$O(1)$ loop model, generalizing the conjectures in~\cite{Raz-Str-2}.
An example of a domain in one of these special subfamilies is given in
figure~\ref{fig.multicone}, right.

\appendix
\section{Glossary of states}
\label{app.gloss}

\noindent
Here we collect pictures describing all the ``states'' (in the linear
space $\mathbb{C}^{\LP(n)}$) which are used in the paper. This is
intended as a glossary, collecting all the definitions scattered
within the text.

\begin{figure}[ht!]
\[
\ket{s}\quad
\setlength{\unitlength}{10pt}
\raisebox{-20pt}{
\begin{picture}(14,4.5)(-.5,0)
\put(0.25,0.5){\includegraphics[scale=2]{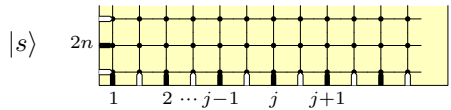}}
\put(-0.2,2.1){\makebox[0pt][c]{$\scriptstyle{2n}$}}
\put(1,0){\makebox[0pt][c]{$\scriptstyle{1}$}}
\put(3,0){\makebox[0pt][c]{$\scriptstyle{2}$}}
\put(3.8,0){\makebox[0pt][c]{$\scriptstyle{\cdots}$}}
\put(5,0){\makebox[0pt][c]{$\scriptstyle{j-1}$}}
\put(7,0){\makebox[0pt][c]{$\scriptstyle{j}$}}
\put(9,0){\makebox[0pt][c]{$\scriptstyle{j+1}$}}
\end{picture}
}
\]
\caption{The state $\ket{s}$, collecting all the FPL in a $n \times n$
  square. Here only the bottom part of the domain is shown, and the
  righmost part is left undetermined, in order to treat in an unitary
  way the case of even and odd~$n$.\label{fig.stateS}}
\end{figure}

\begin{figure}[ht!]
\[
\begin{array}{cc}
\ssb{a}{j} 
&
\ssc{a}{j} 
\\
\setlength{\unitlength}{10pt}
\begin{picture}(14,4.5)(-.5,0)
\put(0.25,0.5){\includegraphics[scale=2]{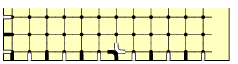}}
\put(-0.2,2.1){\makebox[0pt][c]{$\scriptstyle{2n}$}}
\put(1,0){\makebox[0pt][c]{$\scriptstyle{1}$}}
\put(3,0){\makebox[0pt][c]{$\scriptstyle{2}$}}
\put(3.8,0){\makebox[0pt][c]{$\scriptstyle{\cdots}$}}
\put(5,0){\makebox[0pt][c]{$\scriptstyle{j-1}$}}
\put(7,0){\makebox[0pt][c]{$\scriptstyle{j}$}}
\put(9,0){\makebox[0pt][c]{$\scriptstyle{j+1}$}}
\end{picture}
&
\setlength{\unitlength}{10pt}
\begin{picture}(14,4.5)(-.5,0)
\put(0.25,0.5){\includegraphics[scale=2]{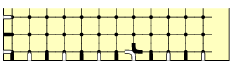}}
\put(-0.2,2.1){\makebox[0pt][c]{$\scriptstyle{2n}$}}
\put(1,0){\makebox[0pt][c]{$\scriptstyle{1}$}}
\put(3,0){\makebox[0pt][c]{$\scriptstyle{2}$}}
\put(3.8,0){\makebox[0pt][c]{$\scriptstyle{\cdots}$}}
\put(5,0){\makebox[0pt][c]{$\scriptstyle{j-1}$}}
\put(7,0){\makebox[0pt][c]{$\scriptstyle{j}$}}
\put(9,0){\makebox[0pt][c]{$\scriptstyle{j+1}$}}
\put(11,0){\makebox[0pt][c]{$\scriptstyle{j+2}$}}
\end{picture}
\\
\setlength{\unitlength}{10pt}
\begin{picture}(14,4.5)(-.5,0)
\put(0.25,0.5){\includegraphics[scale=2]{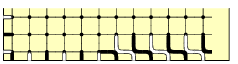}}
\put(-0.2,2.1){\makebox[0pt][c]{$\scriptstyle{2n}$}}
\put(1,0){\makebox[0pt][c]{$\scriptstyle{1}$}}
\put(3,0){\makebox[0pt][c]{$\scriptstyle{2}$}}
\put(3.8,0){\makebox[0pt][c]{$\scriptstyle{\cdots}$}}
\put(5,0){\makebox[0pt][c]{$\scriptstyle{j-1}$}}
\put(7,0){\makebox[0pt][c]{$\scriptstyle{j}$}}
\put(9,0){\makebox[0pt][c]{$\scriptstyle{j+1}$}}
\end{picture}
&
\setlength{\unitlength}{10pt}
\begin{picture}(14,4.5)(-.5,0)
\put(0.25,0.5){\includegraphics[scale=2]{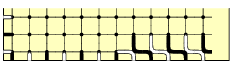}}
\put(-0.2,2.1){\makebox[0pt][c]{$\scriptstyle{2n}$}}
\put(1,0){\makebox[0pt][c]{$\scriptstyle{1}$}}
\put(3,0){\makebox[0pt][c]{$\scriptstyle{2}$}}
\put(3.8,0){\makebox[0pt][c]{$\scriptstyle{\cdots}$}}
\put(5,0){\makebox[0pt][c]{$\scriptstyle{j-1}$}}
\put(7,0){\makebox[0pt][c]{$\scriptstyle{j}$}}
\put(9,0){\makebox[0pt][c]{$\scriptstyle{j+1}$}}
\put(11,0){\makebox[0pt][c]{$\scriptstyle{j+2}$}}
\end{picture}
\\
\setlength{\unitlength}{10pt}
\begin{picture}(14,4.5)(-.5,0)
\put(0.25,0.5){\includegraphics[scale=2]{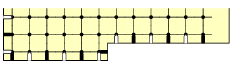}}
\put(-0.2,2.1){\makebox[0pt][c]{$\scriptstyle{2n}$}}
\put(1,0){\makebox[0pt][c]{$\scriptstyle{1}$}}
\put(3,0){\makebox[0pt][c]{$\scriptstyle{2}$}}
\put(3.8,0){\makebox[0pt][c]{$\scriptstyle{\cdots}$}}
\put(5,0){\makebox[0pt][c]{$\scriptstyle{j-1}$}}
\put(6.8,0.75){\makebox[0pt][c]{$\scriptstyle{j}$}}
\put(8.2,1){\makebox[0pt][c]{$\scriptstyle{j+1}$}}
\end{picture}
&
\setlength{\unitlength}{10pt}
\begin{picture}(14,4.5)(-.5,0)
\put(0.25,0.5){\includegraphics[scale=2]{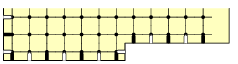}}
\put(-0.2,2.1){\makebox[0pt][c]{$\scriptstyle{2n}$}}
\put(1,0){\makebox[0pt][c]{$\scriptstyle{1}$}}
\put(3,0){\makebox[0pt][c]{$\scriptstyle{2}$}}
\put(3.8,0){\makebox[0pt][c]{$\scriptstyle{\cdots}$}}
\put(5,0){\makebox[0pt][c]{$\scriptstyle{j-1}$}}
\put(7,0){\makebox[0pt][c]{$\scriptstyle{j}$}}
\put(8.4,1){\makebox[0pt][c]{$\scriptstyle{j+1}$}}
\put(10.2,1){\makebox[0pt][c]{$\scriptstyle{j+2}$}}
\end{picture}
\end{array}
\]
\caption{The states $\ssb{a}{j}$ and $\ssc{a}{j}$, here represented
  three times: top, as their basic definition in the $n \times n$
  square, with a constrained site; middle: showing the sites in the
  square that are frozen by the constraint; bottom: the resulting
  geometry obtained removing the frozen part.\label{fig.stateSa}}
\end{figure}

\begin{figure}[ht!]
\[
\begin{array}{cc}
\ssb{b}{j} 
&
\ssc{b}{j} 
\\
\setlength{\unitlength}{10pt}
\begin{picture}(14,4.5)(-.5,0)
\put(0.25,0.5){\includegraphics[scale=2]{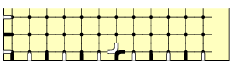}}
\put(-0.2,2.1){\makebox[0pt][c]{$\scriptstyle{2n}$}}
\put(1,0){\makebox[0pt][c]{$\scriptstyle{1}$}}
\put(3,0){\makebox[0pt][c]{$\scriptstyle{2}$}}
\put(3.8,0){\makebox[0pt][c]{$\scriptstyle{\cdots}$}}
\put(5,0){\makebox[0pt][c]{$\scriptstyle{j-1}$}}
\put(7,0){\makebox[0pt][c]{$\scriptstyle{j}$}}
\put(9,0){\makebox[0pt][c]{$\scriptstyle{j+1}$}}
\end{picture}
&
\setlength{\unitlength}{10pt}
\begin{picture}(14,4.5)(-.5,0)
\put(0.25,0.5){\includegraphics[scale=2]{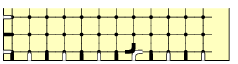}}
\put(-0.2,2.1){\makebox[0pt][c]{$\scriptstyle{2n}$}}
\put(1,0){\makebox[0pt][c]{$\scriptstyle{1}$}}
\put(3,0){\makebox[0pt][c]{$\scriptstyle{2}$}}
\put(3.8,0){\makebox[0pt][c]{$\scriptstyle{\cdots}$}}
\put(5,0){\makebox[0pt][c]{$\scriptstyle{j-1}$}}
\put(7,0){\makebox[0pt][c]{$\scriptstyle{j}$}}
\put(9,0){\makebox[0pt][c]{$\scriptstyle{j+1}$}}
\put(11,0){\makebox[0pt][c]{$\scriptstyle{j+2}$}}
\end{picture}
\\
\setlength{\unitlength}{10pt}
\begin{picture}(14,4.5)(-.5,0)
\put(0.25,0.5){\includegraphics[scale=2]{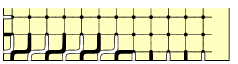}}
\put(-0.2,2.1){\makebox[0pt][c]{$\scriptstyle{2n}$}}
\put(1,0){\makebox[0pt][c]{$\scriptstyle{1}$}}
\put(3,0){\makebox[0pt][c]{$\scriptstyle{2}$}}
\put(3.8,0){\makebox[0pt][c]{$\scriptstyle{\cdots}$}}
\put(5,0){\makebox[0pt][c]{$\scriptstyle{j-1}$}}
\put(7,0){\makebox[0pt][c]{$\scriptstyle{j}$}}
\put(9,0){\makebox[0pt][c]{$\scriptstyle{j+1}$}}
\end{picture}
&
\setlength{\unitlength}{10pt}
\begin{picture}(14,4.5)(-.5,0)
\put(0.25,0.5){\includegraphics[scale=2]{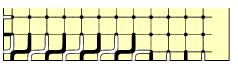}}
\put(-0.2,2.1){\makebox[0pt][c]{$\scriptstyle{2n}$}}
\put(1,0){\makebox[0pt][c]{$\scriptstyle{1}$}}
\put(3,0){\makebox[0pt][c]{$\scriptstyle{2}$}}
\put(3.8,0){\makebox[0pt][c]{$\scriptstyle{\cdots}$}}
\put(5,0){\makebox[0pt][c]{$\scriptstyle{j-1}$}}
\put(7,0){\makebox[0pt][c]{$\scriptstyle{j}$}}
\put(9,0){\makebox[0pt][c]{$\scriptstyle{j+1}$}}
\put(11,0){\makebox[0pt][c]{$\scriptstyle{j+2}$}}
\end{picture}
\\
\setlength{\unitlength}{10pt}
\begin{picture}(14,4.5)(-.5,0)
\put(0.25,0.5){\includegraphics[scale=2]{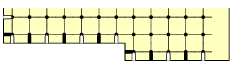}}
\put(-0.2,2.1){\makebox[0pt][c]{$\scriptstyle{2n}$}}
\put(2,1){\makebox[0pt][c]{$\scriptstyle{1}$}}
\put(4,1){\makebox[0pt][c]{$\scriptstyle{2}$}}
\put(4.8,1){\makebox[0pt][c]{$\scriptstyle{\cdots}$}}
\put(6,1){\makebox[0pt][c]{$\scriptstyle{j-1}$}}
\put(7.1,0.75){\makebox[0pt][c]{$\scriptstyle{j}$}}
\put(9,0){\makebox[0pt][c]{$\scriptstyle{j+1}$}}
\put(10,0){\makebox[0pt][c]{$\scriptstyle{\cdots}$}}
\end{picture}
&
\setlength{\unitlength}{10pt}
\begin{picture}(14,4.5)(-.5,0)
\put(0.25,0.5){\includegraphics[scale=2]{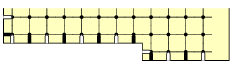}}
\put(-0.2,2.1){\makebox[0pt][c]{$\scriptstyle{2n}$}}
\put(2,1){\makebox[0pt][c]{$\scriptstyle{1}$}}
\put(4,1){\makebox[0pt][c]{$\scriptstyle{2}$}}
\put(4.8,1){\makebox[0pt][c]{$\scriptstyle{\cdots}$}}
\put(6,1){\makebox[0pt][c]{$\scriptstyle{j-1}$}}
\put(8,1){\makebox[0pt][c]{$\scriptstyle{j}$}}
\put(9,0){\makebox[0pt][c]{$\scriptstyle{j+1}$}}
\put(11,0){\makebox[0pt][c]{$\scriptstyle{j+2}$}}
\put(12,0){\makebox[0pt][c]{$\scriptstyle{\cdots}$}}
\end{picture}
\end{array}
\]
\caption{The states $\ssb{b}{j}$ and $\ssc{b}{j}$, here represented
  three times, with the same notations as in
  figure~\ref{fig.stateSa}.\label{fig.stateSb}}
\end{figure}

\begin{figure}[ht!]
\[
\begin{array}{cc}
\ssb{c}{j} 
&
\ssc{c}{j} 
\\
\setlength{\unitlength}{10pt}
\begin{picture}(14,4.5)(-.5,0)
\put(0.25,0.5){\includegraphics[scale=2]{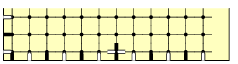}}
\put(-0.2,2.1){\makebox[0pt][c]{$\scriptstyle{2n}$}}
\put(1,0){\makebox[0pt][c]{$\scriptstyle{1}$}}
\put(3,0){\makebox[0pt][c]{$\scriptstyle{2}$}}
\put(3.8,0){\makebox[0pt][c]{$\scriptstyle{\cdots}$}}
\put(5,0){\makebox[0pt][c]{$\scriptstyle{j-1}$}}
\put(7,0){\makebox[0pt][c]{$\scriptstyle{j}$}}
\put(9,0){\makebox[0pt][c]{$\scriptstyle{j+1}$}}
\end{picture}
&
\setlength{\unitlength}{10pt}
\begin{picture}(14,4.5)(-.5,0)
\put(0.25,0.5){\includegraphics[scale=2]{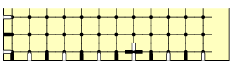}}
\put(-0.2,2.1){\makebox[0pt][c]{$\scriptstyle{2n}$}}
\put(1,0){\makebox[0pt][c]{$\scriptstyle{1}$}}
\put(3,0){\makebox[0pt][c]{$\scriptstyle{2}$}}
\put(3.8,0){\makebox[0pt][c]{$\scriptstyle{\cdots}$}}
\put(5,0){\makebox[0pt][c]{$\scriptstyle{j-1}$}}
\put(7,0){\makebox[0pt][c]{$\scriptstyle{j}$}}
\put(9,0){\makebox[0pt][c]{$\scriptstyle{j+1}$}}
\put(11,0){\makebox[0pt][c]{$\scriptstyle{j+2}$}}
\end{picture}
\\
\setlength{\unitlength}{10pt}
\begin{picture}(14,4.5)(-.5,0)
\put(0.25,0.5){\includegraphics[scale=2]{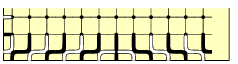}}
\put(-0.2,2.1){\makebox[0pt][c]{$\scriptstyle{2n}$}}
\put(1,0){\makebox[0pt][c]{$\scriptstyle{1}$}}
\put(3,0){\makebox[0pt][c]{$\scriptstyle{2}$}}
\put(3.8,0){\makebox[0pt][c]{$\scriptstyle{\cdots}$}}
\put(5,0){\makebox[0pt][c]{$\scriptstyle{j-1}$}}
\put(7,0){\makebox[0pt][c]{$\scriptstyle{j}$}}
\put(9,0){\makebox[0pt][c]{$\scriptstyle{j+1}$}}
\end{picture}
&
\setlength{\unitlength}{10pt}
\begin{picture}(14,4.5)(-.5,0)
\put(0.25,0.5){\includegraphics[scale=2]{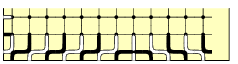}}
\put(-0.2,2.1){\makebox[0pt][c]{$\scriptstyle{2n}$}}
\put(1,0){\makebox[0pt][c]{$\scriptstyle{1}$}}
\put(3,0){\makebox[0pt][c]{$\scriptstyle{2}$}}
\put(3.8,0){\makebox[0pt][c]{$\scriptstyle{\cdots}$}}
\put(5,0){\makebox[0pt][c]{$\scriptstyle{j-1}$}}
\put(7,0){\makebox[0pt][c]{$\scriptstyle{j}$}}
\put(9,0){\makebox[0pt][c]{$\scriptstyle{j+1}$}}
\put(11,0){\makebox[0pt][c]{$\scriptstyle{j+2}$}}
\end{picture}
\\
\setlength{\unitlength}{10pt}
\begin{picture}(14,4.5)(-.5,0)
\put(0.25,0.5){\includegraphics[scale=2]{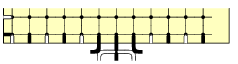}}
\put(-0.2,2.1){\makebox[0pt][c]{$\scriptstyle{2n}$}}
\put(2,1){\makebox[0pt][c]{$\scriptstyle{1}$}}
\put(4,1){\makebox[0pt][c]{$\scriptstyle{2}$}}
\put(5,1){\makebox[0pt][c]{$\scriptstyle{\cdots}$}}
\put(5.2,0.3){\makebox[0pt][c]{$\scriptstyle{j-1}$}}
\put(7,0){\makebox[0pt][c]{$\scriptstyle{j}$}}
\put(8.8,0.3){\makebox[0pt][c]{$\scriptstyle{j+1}$}}
\put(10,1){\makebox[0pt][c]{$\scriptstyle{j+2}$}}
\put(11,1){\makebox[0pt][c]{$\scriptstyle{\cdots}$}}
\end{picture}
&
\setlength{\unitlength}{10pt}
\begin{picture}(14,4.5)(-.5,0)
\put(0.25,0.5){\includegraphics[scale=2]{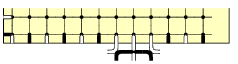}}
\put(-0.2,2.1){\makebox[0pt][c]{$\scriptstyle{2n}$}}
\put(2,1){\makebox[0pt][c]{$\scriptstyle{1}$}}
\put(4,1){\makebox[0pt][c]{$\scriptstyle{2}$}}
\put(4.7,1){\makebox[0pt][c]{$\scriptstyle{\cdots}$}}
\put(5.7,1){\makebox[0pt][c]{$\scriptstyle{j-1}$}}
\put(7,0){\makebox[0pt][c]{$\scriptstyle{j}$}}
\put(9,0){\makebox[0pt][c]{$\scriptstyle{j+1}$}}
\put(10.3,1){\makebox[0pt][c]{$\scriptstyle{j+2}$}}
\put(11.5,1){\makebox[0pt][c]{$\scriptstyle{\cdots}$}}
\end{picture}
\end{array}
\]
\caption{The states $\ssb{c}{j}$ and $\ssc{c}{j}$, here represented
  three times, with the same notations as in
  figure~\ref{fig.stateSa}.\label{fig.stateSc}}
\end{figure}

\section*{Acknowledgements}

\noindent
We thank P.~Di Francesco, P.~Zinn-Justin and J.-B.~Zuber for a
critical analysis of the proof.
L.C.~also thanks his wife G\'eraldine for her patience during the
turbulent writing of this manuscript.

L.C.~acknowledges the financial support of the ANR program
``SLE'', ANR-06-BLAN-0058-01.


\end{document}